# Large time behaviors of upwind schemes by jump processes

## Lei Li[*1] and Jian-Guo Liu[†2]

[1]Institute of Natural Sciences, Shanghai Jiao Tong University, Minhang District, Shanghai 200240, China.
[2] Department of Mathematics and Department of Physics, Duke University, Durham, NC 27708, USA.

### Abstract

We revisit the traditional upwind schemes for linear conservation laws in the view-point of jump processes, allowing studying upwind schemes using probabilistic tools. In particular, for Fokker-Planck equations on $\mathbb{R}$, in the case of weak confinement, we show that the solution of upwind scheme converges to a stationary solution. In the case of strong confinement, using a discrete Poincaré inequality, we prove that the $O(h)$ numeric error under $\ell^1$ norm is uniform in time, and establish the uniform exponential convergence to the steady states. Compared with the traditional results of exponential convergence of upwind schemes, our result is in the whole space without boundary. We also establish similar results on torus for which the stationary solution of the scheme does not have detailed balance. This work shows an interesting connection between standard numerical methods and time continuous Markov chains, and could motivate better understanding of numerical analysis for conservation laws.

## 1 Introduction

In discretizing hyperbolic equations or the convection terms in mixed type equations such as Navier-Stokes equations and conservation laws with diffusion, the upwind scheme is usually used to numerically simulate the direction of propagation of information to ensure desired stability [1]. For nonlinear hyperbolic conservation laws, the upwind schemes are especially important. In fact, for scalar conservation laws, when the time and spatial discretization satisfies the Courant–Friedrichs–Lewy (CFL) condition, the upwind schemes lead to the so-called monotone schemes [2]. It is well known that monotone schemes guarantee that the numerical solutions converge to the entropy weak solution [2, 3]. The entropy weak solution is important for physical phenomena like shocks. Monotone schemes are often first order accurate and were later extended to total variation diminishing (TVD) schemes satisfying the entropy inequality [4, 5]. It was proved that TVD schemes satisfying the entropy inequality can yield numerical solutions converging to the unique entropy weak solution on $[0, T] \times D$, where $D$ is a compact set in space [4, 5]. The proofs of existence of the weak solutions in [3, 4, 5] relies on the boundedness of variation and $L^1$ norms, which imply compactness in $L^1([0, T] \times D)$. Since TVD guarantees the compactness in $L^1$, it is an important property for schemes of conservation laws. One can consider method-of-lines form (with time being continuous) that are TVD and then find time discretization that preserves TVD [6]. As we will see soon, method-of-line scheme with upwind differencing in space is TVD.

---

[*]leiliyushan@gmail.com
[†]jliu@phy.duke.edu

Now, let us focus this in the 1D case. In general, the scalar conservation law in 1D space is given by

$$\partial_t \rho + f(x, \rho)_x = \partial_x(D(x)\partial_x \rho). \tag{1.1}$$

We will assume all functions are smooth enough, $f(x, 0) = 0$ and $D(x) \geq 0$. If $D(x) = 0$, we have the standard hyperbolic conservation laws. For $D = 0$, Kružkov proved in [7] that if $\partial_x f(x, \rho)$ is locally Lipschitz in $\rho$, the bounded weak solution satisfying an entropy condition ([7, Definition 1]) is unique. The existence result of such solutions in [7] requires that the derivatives of $f(x, \rho)$ satisfy some boundedness conditions uniform in $x$ so that the vanishing viscosity method works. In particular, if $f(x, \rho) = f_1(\rho)$ with $f_1$ being locally Lipschitz, the existence result holds. With suitable assumptions on the flux $f(x, \rho)$, like $f(x, \rho) = f_1(\rho)$, or some confinement conditions, $\int_{\mathbb{R}} \rho \, dx$ is a constant (see, for example, [8, Proposition 2.3.6]). For general fluxes that can depend on $x$, even if the equation is well-posed, the total mass can decay because some mass can escape to infinity, like $\rho_t + \partial_x((1 + x^2)\rho) = 0$.

In this paper, we look at the upwind schemes in the method-of-lines form. For simulation, one can use the methods in [6] to get fully discretized TVD schemes or just leave the time variable continuous as in section 7 . For upwind discretization, we decompose the flux as

$$f = f_+ - f_-, \ \partial_\rho f_\pm(x, \rho) \geq 0, \ f_\pm(x, 0) = 0, \ i = 1, 2. \tag{1.2}$$

Clearly, we can set

$$f_\pm(x, \rho) = \int_0^\rho (\partial_\rho f(x, v))^\pm \, dv, \tag{1.3}$$

where we have used $z^+ = z \vee 0$ and $z^- = -z \wedge 0$ for $z \in \mathbb{R}$. If $f \in C^1$, $f_\pm$ is also $C^1$.

We discretize the space with step size $h > 0$ and set $x_j = jh$. Let $\rho_j(t)$ be the numerical solution at site $x_j$, with $\rho_j(0)$ being some approximation for $\frac{1}{h} \int_{x_{j-1/2}}^{x_{j+1/2}} \rho(x, 0) \, dx$. Then, the upwind scheme for (1.1) can be constructed based on the flux splitting [3, 9]

$$\frac{d}{dt}\rho_j = -\left(\frac{f^+(x_j, \rho_j) - f^+(x_{j-1}, \rho_{j-1})}{h} - \frac{f^-(x_{j+1}, \rho_{j+1}) - f^-(x_j, \rho_j)}{h}\right)$$
$$+ \frac{1}{h^2}(D_{j+1/2}\rho_{j+1} - (D_{j+1/2} + D_{j-1/2})\rho_j + D_{j-1/2}\rho_{j-1}), \tag{1.4}$$

where $D_{j+1/2} = D(x_j + \frac{h}{2})$. We denote $f_{\pm, j} := f_\pm(x_j, \rho_j)$. The time continuous upwind scheme (1.4) is TVD for bounded $\ell^1$ solutions that decay fast enough. In other words, if $\rho \in L^\infty(0, T; \ell^1 \cap \ell^\infty)$ is a solution that decays fast enough, $\sum_j |\rho_{j+1} - \rho_j|$ is non-increasing. Here, $L^\infty(0, T; X)$ means the $\|\cdot\|_X$ norm is essentially bounded on $[0, T]$ while $\ell^p$ refers to the usual Banach spaces in numerical analysis (note that there is $h$ involved)

$$\ell^p := \begin{cases} \{\rho : \mathbb{Z} \to \mathbb{R} \ \Big| \ \|\rho\|_{\ell^p} := (\sum_{j \in \mathbb{Z}} h|\rho_j|^p)^{1/p} < \infty\}, \ p \in [1, \infty), \\ \{\rho : \mathbb{Z} \to \mathbb{R} \ \Big| \ \|\rho\|_{\ell^\infty} := \sup_{j \in \mathbb{Z}} |\rho_j| < \infty\}, \ p = \infty. \end{cases} \tag{1.5}$$

The reason that the scheme is TVD is that the numbers

$$a_j^+ := \frac{f^+(x_j, \rho_j) - f^+(x_{j-1}, \rho_{j-1})}{\rho_j - \rho_{j-1}}, \ a_j^- := \frac{f^-(x_{j+1}, \rho_{j+1}) - f^-(x_j, \rho_j)}{\rho_{j+1} - \rho_j}$$

are bounded for given $j$ (since we have assumed $\rho$ is bounded) and non-negative. Using similar technique as in the proof of Proposition 4.1, we can conclude the TVD property. If the TVD property is satisfied, then we can get the convergence of the numerical scheme by compactness in $L^1_{\text{loc}}(\mathbb{R})$. Of course, whether the true solutions of (1.4) decay fast enough depends on concrete conditions on $f(x, \rho)$ and $D(x)$. The upwind scheme (1.4) can be rearranged to the conservative scheme

$$\frac{d}{dt}\rho_j + \frac{1}{h}[J_{j+1/2} - J_{j-1/2}] = 0 \tag{1.6}$$

where

$$J_{j+1/2} = h\alpha_j\rho_j - h\beta_{j+1}\rho_{j+1}, \tag{1.7}$$

with

$$\alpha_j = \frac{f_{+,j}/\rho_j}{h} + \frac{1}{h^2}D_{j+1/2}, \ \ \beta_j = \frac{f_{-,j}/\rho_j}{h} + \frac{1}{h^2}D_{j-1/2}. \tag{1.8}$$

Hence, it can be further written as the discrete form

$$\frac{d}{dt}\rho_j = \alpha_{j-1}\rho_{j-1} + \beta_{j+1}\rho_{j+1} - (\alpha_j + \beta_j)\rho_j. \tag{1.9}$$

According to (1.3), we have for any $j \in \mathbb{Z}$, $f_{\pm,j}/\rho_j \geq 0$ and is bounded for bounded $\rho_j$. If $\rho_j = 0$, the quotient is understood as the partial derivative of $f_{\pm}$ on $\rho$ at $(x_j, 0)$. Hence, the upwind scheme ensures that $\alpha_j, \beta_j$ are nonnegative. We can then interpret the upwind scheme as the master equation of some transition phenomena. In particular, $\alpha_j$ can be understood as the rate of moving the mass from site $j$ to site $j + 1$ while $\beta_j$ the the rate of moving mass from $j$ to $j - 1$. Then (1.9) describes the evolution of mass. Due to this physical understanding, if the upwind scheme (1.6)-(1.7) is well-posed, we expect that (1.9) is non-negativity preserving, and is $\ell^1$ contracting (i.e. $\|\rho^1(t) - \rho^2(t)\|_{\ell^1} \leq \|\rho^1(0) - \rho^2(0)\|_{\ell^1}$). The $\ell^1$ contraction then will imply the uniqueness of solutions.

As we have seen, the upwind difference in space gives TVD, non-negativity preserving and $\ell^1$ contracting schemes (at least in the formal way since the well-posedness needs further investigation). These properties make upwind difference useful in numerical analysis and will guarantee the convergence of numerical schemes on $\mathbb{R} \times [0, T]$. (Indeed, as we will see in sections 3 and 4, for the problems we consider, these properties hold and the mass is also conserved.)

The convergence on $\mathbb{R} \times [0, T]$, however, is not enough if we care about the asymptotic behaviors. Our observation is that when the equation is linear, the master equation (1.9) can be regarded as the forward equation of a jump process (time continuous random walk) [10]. In this case, we can normalize $\rho$ to the probability measure of the random walk on $\mathbb{Z}$. Since jump processes are well-studied [10, 11] in the community of probability, we may then use tools from probability to study the large time behaviors of the upwind schemes. In fact, if the diffusion coefficient is nonzero, it is then the Fokker-Planck equation of a stochastic differential equation (SDE) and we are able to establish the discrete Poincaré inequality under some assumptions inspired by theories in [10, 11]. Hence, in the remaining part of this paper, we only focus on Fokker-Planck equations and investigate the asymptotic behaviors of the upwind schemes using jump processes. We prove the uniform $O(h)$ error of the upwind scheme and prove the exponential convergence to equilibrium states. We remark that the existing results regarding exponential convergence for discrete schemes of conservation laws are often on finite domains (see [12]). We hope our work will bring more understanding of both Markov chains and numerical schemes.

For related references, one can refer, for example, to [13, 14, 15, 16]. Donsker invariance principle [13, 14] claims that a certain rescaled random walk converges to the standard Brownian motion on time interval [0, 1] in distribution. In [15, 16], Markov chains have been used to approximate diffusion processes and the weak convergence of the scheme on fixed time interval has been proved. The general motivation of our work shares similarity with those in [15, 16].

The rest of the paper is organized as follows. In section 2, we give a brief introduction to stochastic differential equations (SDEs) and the associated Fokker-Planck equation. We also have a review of results regarding the stationary distribution and ergodicity. In section 3, we move on to the upwind schemes for the Fokker-Planck equations on $\mathbb{R}$ and show the uniform error estimates. In section 4, we prove some elementary properties of the jump process for the upwind scheme. In particular, we show some basic properties of the discrete backward equation of the Markov jump process and show that the solution of upwind scheme converges to a stationary solution in the case of weak confinement. In section 5, we focus on the strong confinement and study the asymptotic behaviors of the upwind schemes. We show the uniform geometric convergence to the steady states using a discrete Poincaré inequality on the whole space. We then prove the $O(h)$ accuracy for the stationary solution, proving the unproved claim (Theorem 3.1) in section 3. Further in section 6, we establish the results on torus for which detailed balance may not hold. Last in section 7, we propose a Monte Carlo method to numerically solve the upwind schemes in a probabilistic way.

## 2  Preliminaries: basic facts of SDEs

Above, we have mentioned that the linear conservation law with positive diffusion is the Fokker-Planck equation for an SDE. This is the focus of this paper, so we will have a brief review of SDEs in this section. We will focus on general dimension $d$ in this section.

### 2.1  Basic setup of SDEs

The time homogeneous SDE driven by Wiener process in Itô sense is given by [17]:

$$dX = b(X)\,dt + \sigma(X)\,dW. \tag{2.1}$$

The functions $b$ and $\sigma$ are called the drift and diffusion coefficients respectively. $W$ is the standard Wiener process defined on some probability space $(\Omega, \mathcal{F}, \mathbb{P})$. When $b$ and $\sigma$ are Lipschitz continuous and have linear growth at infinity, (2.1) has global strong solutions [17, sect. 5.2] for $L^2(\mathbb{P})$ initial data. The conditions imposed $b(\cdot)$ in [17, sect. 5.2] is too strong for many applications. In fact, it is also known that locally Lipschitz and confinement conditions can also imply the existence and uniqueness of solutions (For example, in [18, Theorem 2.3.5], it is shown that $\max(x \cdot b(x), |\sigma|^2) \le C_1 + C_2 |x|^2$ is enough for the well-posedness, which allows $b$ like $-(1 + |x|^2)^p x$).

The most frequently used confinement condition in this work is the following.

*Assumption* 2.1. Suppose $b$ and $\sigma$ are smooth. The function $b$ satisfies

$$b(x) \cdot x \le -r|x|^2 \tag{2.2}$$

when $|x| > R$ for some $R$. Also, $\sigma$ satisfies $\|\sigma\|_\infty < \infty$ and $\Lambda = \sigma\sigma^T \ge S_1 I > 0$.

Besides this, we sometimes weaken the conditions as follows.

*Assumption* 2.2. Suppose $b$ and $\sigma$ are smooth. The function $b$ satisfies

$$\lim_{|x| \to \infty} \frac{-b(x) \cdot x}{|x|} = \infty. \tag{2.3}$$

Also, $\sigma$ satisfies $\|\sigma\|_\infty < \infty$ and $\Lambda = \sigma\sigma^T \ge S_1 I > 0$.

We will use $\mathbb{E}$ to represent the expectation under $\mathbb{P}$. The notation $\mathbb{E}_x$ indicates that the expectation is conditioned on $X(0) = x$. Let $\mu_t$ be the law of $X(t)$, which is a measure in $\mathbb{R}^d$. Then we have

$$\mathbb{E}f(X(t)) = \langle \mu_t, f \rangle = \int_{\mathbb{R}^d} f\,d\mu_t. \tag{2.4}$$

For smooth bounded function $f(x)$, define

$$u(x, t) = \mathbb{E}_x f(X_t). \tag{2.5}$$

By Itô's calculus [17], $u$ satisfies

$$\partial_t u(x, t) = \mathbb{E}_x \mathcal{L} f(X_t), \tag{2.6}$$

where $\mathcal{L}$ is the generator of the process

$$\mathcal{L} := b \cdot \nabla + \frac{1}{2}\Lambda_{ij}\partial_{ij}, \tag{2.7}$$

where we used Einstein summation convention (i.e. $\Lambda_{ij}\partial_{ij} \equiv \sum_{i,j=1}^{d} \Lambda_{ij}\partial_{x_i x_j}$) and

$$\Lambda = \sigma\sigma^T. \tag{2.8}$$

This is a special case of Dynkin's formula. The density of the law of $X(t)$ starting $x$, denoted by $p(t, x, y)$, is called the Green's function. When $\Lambda$ is positive definite, $p(t, x, y)$ is a smooth function for $t > 0$. Equation (2.6) implies that $p(t, x, y)$ satisfies the forward Kolmogorov equation, or Fokker-Planck equation for $t > 0$:

$$\partial_t p = -\nabla_y \cdot (b(y)p) + \frac{1}{2}\partial_{y_i y_j}(\Lambda_{ij}(y)p) := \mathcal{L}_y^* p, \tag{2.9}$$

where the subindex $y$ means that the derivatives are taken on $y$ variable. By the well-posedness of (2.1), we have under the confinement conditions that

$$\int_{\mathbb{R}^d} p(t, x, y)\, dy = 1, \ \forall x \in \mathbb{R}^d, t > 0. \tag{2.10}$$

Clearly, for general starting probability measure $\mu_0$, the law of $X(t)$ also satisfies (2.9) in the distributional sense:

$$\frac{d}{dt}\langle \mu_t, f \rangle = \langle \mu_t, \mathcal{L}f \rangle,$$

which is clearly a generalization of (2.6). Moreover, let $v : (x, t) \mapsto v(x, t)$ solve the backward Kolmogorov equation

$$\partial_t v = \mathcal{L}v = b \cdot \nabla v + \frac{1}{2}\Lambda_{ij}\partial_{ij}v \tag{2.11}$$

with initial condition $v(x, 0) = f(x)$. Let $X(t)$ be the process satisfying (2.1) with initial condition $X(0) = x$. We check that $Y_s = v(X(s), t - s)$ is a martingale and therefore

$$v(x, t) = Y_0 = \mathbb{E}Y_t = \mathbb{E}v(X(t), 0) = \mathbb{E}f(X(t)) = u(x, t). \tag{2.12}$$

This means that (2.5) solves the backward Kolmogorov equation. Combining with (2.6), we can infer that the Green's function satisfies $\mathcal{L}_y^* p(t, x, y) = \mathcal{L}_x p(t, x, y)$, or

$$-\nabla_y \cdot (b(y)p(t, x, y)) + \frac{1}{2}\partial_{y_i y_j}(\Lambda_{ij}(y)p(t, x, y)) =$$
$$b(x) \cdot \nabla_x p(t, x, y) + \frac{1}{2}\Lambda_{ij}(x)\partial_{x_i x_j}p(t, x, y). \tag{2.13}$$

## 2.2 Stationary solutions and ergodicity

Under Assumption 2.1, using Itô's formula and test function $f(x) = \exp(c|x|^2)$, one can show that

$$\mathbb{E}_x \exp(c|X_t|^2) \le \exp(c|x|^2)e^{-rt} + C, \tag{2.14}$$

for some positive constants $c, r, C$. This implies that the process has certain recurrent properties so that the SDE (2.1) has a unique stationary distribution $\pi$ [19, sect. 4.4-4.7]. Moreover, $\pi$ has a density with respect to Lebesgue measure [19, Lemma 4.16]. Below, we may abuse the notation a little bit and use $\pi(\cdot)$ to mean this density for convenience. The Green's function $p(t, x, y)$ converges to $\pi(y)$ pointwise as $t \to \infty$ for all $x \in \mathbb{R}^d$ [19, Lemma 4.17]. Clearly, $\pi(y)$ has finite moment of any order by (2.14). Since $\pi(y)$ is a solution to the parabolic equation (2.9) with the diffusion coefficient matrix positive definite, $\pi(y)$ is smooth and $\pi(y) > 0$.

Often people study the ergodicity of SDEs in the $L^p$ spaces. We will use $L^p(\mathbb{R}^d)$ to represent the $L^p$ spaces associated with the Lebesgue measure while $L^p(\nu)$ to mean the $L^p$ spaces associated with the measure $\nu$. If $\nu$ has a density $w$, we also write $L^p(\nu)$ as $L^p(w)$. The most frequently used weight is $w = \pi$. Let $p(\cdot, t)$ be the density of $\mu_t$. We often define

$$q(x, t) := \frac{p(x, t)}{\pi(x)} \geq 0,$$

and study the convergence of $q(\cdot, t)$ to 1 in $L^p(\pi)$ spaces.

Note that $\Lambda_{ij}$ is symmetric and

$$-\nabla \cdot (b\pi) + \frac{1}{2}\partial_{ij}(\Lambda_{ij}\pi) = 0, \tag{2.15}$$

we have

$$\partial_t q = \left(\frac{1}{\pi}\nabla \cdot (\Lambda\pi) - b\right) \cdot \nabla q + \frac{1}{2}\Lambda_{ij}\partial_{ij}q. \tag{2.16}$$

If the detailed balance condition

$$b = \frac{1}{2\pi}\nabla \cdot (\Lambda\pi) \tag{2.17}$$

holds (for example, $\Lambda = 2DI$ and $b = -\nabla V$), which clearly indicates (2.15), then we have the useful identity

$$\mathcal{L}^*(f\pi) = \pi\mathcal{L}f + f\mathcal{L}^*\pi = \pi\mathcal{L}f. \tag{2.18}$$

Then (2.16) can be rewritten as

$$\partial_t q = b \cdot \nabla q + \frac{1}{2}\Lambda_{ij}\partial_{ij}q, \tag{2.19}$$

which is the backward equation (2.11). In this case, the semigroup $e^{t\mathcal{L}}$ is symmetric in $L^2(\pi)$ and $e^{t\mathcal{L}^*}$ is symmetric in $L^2(1/\pi)$ by (2.18). Hence, it is convenient to investigate $u(\cdot, t) \to \langle \pi, f \rangle$ and $q(\cdot, t) \to 1$ in $L^2(\pi)$ using (2.11). If the detailed balance is not satisfied, the modified generator

$$\tilde{\mathcal{L}} = \left(\frac{1}{\pi}\nabla \cdot (\Lambda\pi) - b\right) \cdot \nabla + \frac{1}{2}\Lambda_{ij}\partial_{ij} =: \tilde{b} \cdot \nabla + \frac{1}{2}\Lambda_{ij}\partial_{ij} \tag{2.20}$$

corresponds to another SDE

$$dY = \tilde{b}\,dt + \sigma dY, \tag{2.21}$$

which has the same stationary distribution $\pi$, or $\tilde{\mathcal{L}}^*\pi = 0$. Suppose the law of $X(0)$ has a density $p^0(y)$. It follows from (2.21) that

$$q(x, t) = \mathbb{E}\left(\frac{p^0(Y(t))}{\pi(Y(t))}\Big| Y(0) = x\right). \tag{2.22}$$

Hence, though the semigroups generated by $\mathcal{L}$ and $\tilde{\mathcal{L}}$ are not symmetric in $L^2(\pi)$, one can still consider the convergence of $u(\cdot, t) \to \langle \pi, f \rangle$ and $q(\cdot, t) \to 1$ in $L^2(\pi)$ using Kolmogorov backward equations.

It is well-known that Condition 2.1 implies geometric ergodicity regarding the convergence of $u(\cdot, t)$ to $\langle \pi, f \rangle$ or $\mu_t$ to $\pi$ using coupling argument for SDEs. In particular, we have the $V$-uniform geometric ergodicity for $u(\cdot, t) \to \langle \pi, f \rangle$ ([20, 21]) or geometric convergence of $\mu_t \to \pi$ in Wasserstein space ([22, 23]). Besides the coupling argument, one may prove the geometric convergence of $u(\cdot, t)$ to $\langle \pi, f \rangle$ in $L^p(\pi)$ spaces using spectral gap and Perron-Frobenius type theorems (see [20, Chap. 20]; [24, 25, 26, 27, 28] for example). The $V$-uniform ergodicity and ergodicity in $L^p(\pi)$ do not necessarily imply each other, unless extra conditions are imposed [20, Chap. 20].

The geometric convergence of $q(\cdot, t)$ to 1 (equivalent to the convergence of $u(\cdot, t) \to \langle \pi, f \rangle$ for the modified SDE (2.21)) in $L^p(\pi)$ for the modified SDE (2.21)) in $L^p(\pi)$ spaces can also be obtained directly using the Fokker-Planck equation and some functional inequalities (Poincaré inequality, or log Sobolev inequality etc). These functional inequalities will imply spectral gaps of the semigroups. Let us explain this briefly. Take a smooth function $\varphi$ and recall (2.16). We find

$$\frac{d}{dt} \varphi(q) = \tilde{\mathcal{L}}(\varphi(q)) - \frac{1}{2} \varphi''(q) \Lambda_{ij} \partial_i q \partial_j q.$$

Multiplying $\pi$ and taking integral (recall $\tilde{\mathcal{L}}^*(\pi) = 0$), we have the energy-dissipation relation

$$\frac{d}{dt} \mathcal{F} := \frac{d}{dt} \int_{\mathbb{R}} \varphi(q) \pi \, dx = -\frac{1}{2} \int_{\mathbb{R}} \varphi''(q) \Lambda_{ij} \partial_i q \partial_j q \, \pi dx =: -\mathcal{D}. \qquad (2.23)$$

If $\varphi$ is the quadratic function and the Poincaré inequality associated with $\pi$ can be established, the geometric convergence of $q(\cdot, t)$ to 1 in $L^2(\pi)$ follows. This clearly implies that the geometric convergence of $q(\cdot, t)$ to 1 in $L^1(\pi)$ and hence the geometric convergence of $p(\cdot, t)$ to $\pi$ in $L^1(\mathbb{R}^d)$ norm (total variation norm). Alternatively, one may take $\varphi(q) = q \log q - q + 1$ and then $\mathcal{F}$ becomes the relative entropy or Kullback–Leibler (KL) divergence. If the log-Sobolev inequality holds, one can then establish the geometric convergence of the relative entropy and thus in total variation norm by Pinsker's inequality. The advantage of log-Sobolev inequality is that the constant is dimension free. For the case $b = -\nabla V$ and $\sigma = \sqrt{2D}I$, these results are well-known and one can refer to the review by Markowich and Villani [29].

Below, we focus on the upwind scheme for the Fokker-Planck equation in 1D case only ($d = 1$). The generalization to general $d$ is nontrivial, especially for non-uniform meshes. One issue is that we may lose the detailed balance for discrete schemes and the uniform functional inequalities for general cases are hard to prove. We leave these studies to future.

For $d = 1$, we have the following staightforward observation, which is needed for the error analysis of the upwind scheme:

**Lemma 2.1.** *Let $d = 1$. If $b$ and $\sigma$ satisfy Assumption 2.1, then for any index $n > 0$, there exist positive constants $C_n > 0$, $\nu_n > 0$ such that*

$$\left| \frac{d^n}{dx^n} \pi(x) \right| \leq C_n \exp(-\nu_n |x|^2).$$

To see this, we note that the detailed balance condition $-b\pi + \frac{1}{2} \partial_x(\sigma^2 \pi) = 0$ holds. We can then solve $\sigma^2 \pi$ and therefore $\pi$. Using the formula, Lemma 2.1 follows directly. The details are omitted. For $d > 1$, in the case $b = -\nabla V$ and $\sigma = \sqrt{2D}I$, the claim is also trivial since $\pi \propto \exp(-V/D)$. For general dimension and general $b, \sigma$, we believe Lemma 2.1 is still true due to (2.14) (one may replace the test function $\exp(c|x|^2)$ with the derivatives $x \exp(c|x|^2)$ to get the estimates for derivative of $\pi$). Since we do not need general dimension in this paper, we omit them here. For interested readers, one may refer to [30] for the pointwise estimates at infinity and to, for example, [31, 32, 30] for the theories of elliptic equations in unbounded domains.

# 3 Upwind scheme for Fokker-Planck equation

We can rewrite the Fokker-Planck equation into the conservative form as

$$\partial_t \rho = -\partial_x((b - \sigma \sigma')\rho) + \frac{1}{2} \partial_x(\sigma^2 \partial_x \rho). \qquad (3.1)$$

Here, we assume

$$\rho(x, 0) = \rho^0(x) \geq 0. \qquad (3.2)$$

Clearly, $f(x, \rho) = (b - \sigma \sigma')\rho =: s(x)\rho$. In this case, we have the corresponding decomposition

$$b - \sigma \sigma' =: s^+ - s^-, \quad s^\pm \geq 0. \qquad (3.3)$$

Recall that we use spatial step $h$ to discretize the space and $x_j = jh$. Moreover, we use $R_g : C(\mathbb{R}) \to \mathbb{R}^{\mathbb{Z}}$ to mean the restriction onto the grid:

$$R_g \varphi = (\varphi(x_j)). \tag{3.4}$$

The upwind scheme (1.4) becomes

$$\frac{d}{dt}\rho_j(t) = -\left(\frac{s_j^+ \rho_j - s_{j-1}^+ \rho_{j-1}}{h} - \frac{s_{j+1}^- \rho_{j+1} - s_j^- \rho_j}{h}\right)$$
$$+ \frac{1}{2h^2}(\sigma_{j+1/2}^2 \rho_{j+1} - (\sigma_{j+1/2}^2 + \sigma_{j-1/2}^2)\rho_j + \sigma_{j-1/2}^2 \rho_{j-1}). \tag{3.5}$$

The rates (1.8) for the master equation are independent of $\rho$:

$$\alpha_j = \frac{s_j^+}{h} + \frac{1}{2h^2}\sigma_{j+1/2}^2, \ \ \beta_j = \frac{s_j^-}{h} + \frac{1}{2h^2}\sigma_{j-1/2}^2. \tag{3.6}$$

Denote the sequence

$$\rho^h(t) := (\rho_j(t))_{j\in\mathbb{Z}}. \tag{3.7}$$

We assume $\rho^0(\cdot) \in L^1(\mathbb{R})$ and $\rho_0^h$ is constructed so that

$$\|\rho^h(0) - R_g\rho^0\|_{\ell^1} \le Ch, \ \|\|\rho^h(0)\|_{\ell^1} - \|\rho^0\|_{L^1(\mathbb{R})}\| \le Ch. \tag{3.8}$$

Recall that $\ell^1$ and $L^1(\mathbb{R})$ spaces are introduced in equation (1.5) and section 2.2 respectively. Let $p^0(x) = \rho^0/\|\rho^0\|_{L^1(\mathbb{R})}$. With Assumption 2.2, the SDE (2.1) is not explosive by [18, Theorem 2.3.5]. Hence, $p(x,t)$, the density of the law of $X(t)$, exists and is unique with $\int_{\mathbb{R}} p(x,t)\,dx = 1$. It is the solution of (3.1) with initial condition $p(x,0) = p^0(x)$, and thus $\rho(x,t) = p(x,t)\|\rho^0\|_{L^1(\mathbb{R})}$.

Since the discrete equation is also linear, we can normalize

$$p_j(t) := h\frac{\rho_j(t)}{\|\rho_0^h\|_{\ell^1}} \ge 0 \tag{3.9}$$

so that $p_j(0) \ge 0$ and $\sum_j p_j(0) = 1$. For convenience, we define the sequence

$$p^h(t) := (p_j(t))_{j\in\mathbb{Z}}. \tag{3.10}$$

**Remark 3.1.** Note that $\rho^h(t)$ is the numerical approximation of $\rho(\cdot,t)$, but $p^h(t)$ is not the numerical approximation of the continuous probability density $p(\cdot,t)$ directly. Instead, $h^{-1}p^h(t)$ approximates the probability density $p(\cdot,t)$ and the reason we use this convention shall be clear soon.

The upwind scheme ensures that $\alpha_j, \beta_j$ are nonnegative. Hence, the equation for $p^h(t)$

$$\frac{d}{dt}p_j = \alpha_{j-1}p_{j-1} + \beta_{j+1}p_{j+1} - (\alpha_j + \beta_j)p_j =: (\mathcal{L}_h^* p^h)_j. \tag{3.11}$$

can be regarded as the the forward equation (discrete Fokker-Planck equation) of a jump process or time continuous Markov chain $Z(t)$ [10]. $\alpha_j$ is the rate of jumping from site $j$ to site $j+1$ while $\beta_j$ the the rate of jumping from $j$ to $j-1$.

Here, $\mathcal{L}_h^* : \mathbb{R}^{\mathbb{Z}} \to \mathbb{R}^{\mathbb{Z}}$ is defined for any sequence, but the equation may not have solutions for arbitrarily given initial data. Later in Section 4, we will see that under Assumption 2.2 the chain is nonexplosivethe and equation (3.11) is well-posed for $\ell^1$ initial data. Moreover, $p_j(0) \ge 0$ and $\sum_j p_j(0) = 1$ imply that and that $p_j(t) \ge 0$, $\sum_j p_j(t) = 1$. Then $p_j(t)$ is the probability of appearing at site $j$ and this is why we use the normalization in (3.9).

For the convenience, we define the semigroup as

$$e^{t\mathcal{L}_h^*}p^h(0) := p^h(t). \tag{3.12}$$

With the well-posedness facts and the discussion in the introduction (section 1), we can deduce easily the following, and we omit the proofs.

**Lemma 3.1.** *The semigroup $e^{t\mathcal{L}_h^*}$ for upwind scheme (3.5) or (3.11) is $\ell^1$ contracting and nonnegativity preserving. Moreover, the scheme (3.5) is TVD for $\rho_0^h \in \ell^1$ (i.e. $\sum_j |\rho_j(t) - \rho_{j-1}(t)|$ is non-increasing.)*

In this section, we will focus on the approximation error of the upwind scheme.

## 3.1 Stationary solutions

Consider a stationary solution $\pi^h$ to (3.11) or (1.6). We then find

$$J_{j+1/2}^h = J = const.$$

We take $j \to \infty$ and find $J = 0$. Hence, we have

$$\alpha_j \pi_j^h = \beta_{j+1} \pi_{j+1}^h. \tag{3.13}$$

This is the detailed balance condition.

**Lemma 3.2.** *Suppose the weak confinement 2.2 holds. Then, there is a unique stationary distribution $\pi^h$ for the jump process $Z(t)$ corresponding to (3.11).*

*Proof.* Using (3.13), we find

$$\pi_j^h = \pi_0^h \prod_{k=1}^{j} \frac{\alpha_{k-1}}{\beta_k}$$

With the condition, $b(x_j) < 0$ for $j > 0$ large enough. For these $j$, $\alpha_{j-1}$ is bounded, while $\beta_j$ goes to infinity as $j \to \infty$ by (3.6). Hence, $\pi_j^h$ decays with at least geometric rate. This means $\sum_{j \geq 0} \pi_j^h < \infty$. Similarly, $\sum_{j < 0} \pi_j^h < \infty$ also holds. Hence, $\sum_j \pi_j^h < \infty$ and we can normalize it to a probability distribution so that $\pi_0^h$ is determined uniquely. $\qquad\square$

Similar with $p^h$, $\pi^h$ does not approximate the density $\pi(\cdot)$ of stationary distribution of (2.1). Instead, $h^{-1}\pi^h$ approximates $\pi(\cdot)$. In fact, in section 5.2, we will prove the following, which says that $h^{-1}\pi^h$ approximates $\pi(\cdot)$ with error $h$:

**Theorem 3.1.** *Suppose Assumption 2.1 holds with $S_1 \leq \sigma^2 \leq S_2$. Let $\pi^h$ be the stationary distribution of the jump process $Z(t)$ for (3.11) and $\pi(\cdot)$ be the density of the stationary distribution for (2.1). Then there exist $h_0 > 0$ and $C > 0$ such that (recall (3.4) for $R_g$)*

$$\left\| R_g \pi - \frac{1}{h}\pi^h \right\|_{\ell^1} \leq Ch, \ \forall h \leq h_0. \tag{3.14}$$

On bounded domain, usual techniques for the finite difference method of elliptic equations can be used to prove such type of results. The difference is that now the domain is infinite. The proof relies on the spectral gap of the operator. See Section 5.2 for more details.

Now, let us look at the OU process as an example. We take $b(x) = -x$, $\sigma = 1$. Then,

$$\begin{cases} \alpha_j = \frac{1}{2h^2}, \ \beta_j = j + \frac{1}{2h^2}, \ j \geq 0, \\ \alpha_j = |j| + \frac{1}{2h^2}, \ \beta_j = \frac{1}{2h^2}, \ j < 0. \end{cases}$$

Using the fact $\frac{\pi_{j+1}^h}{\pi_j^h} = \frac{\alpha_j}{\beta_{j+1}}$, we find that $\pi_j^h$ is even. Hence, we only need to focus on $j \geq 0$. Clearly, for $j \geq 1$,

$$\pi_j^h = A_h \prod_{k=1}^{j} \frac{1}{1 + 2kh^2} =: A_h v_j, \tag{3.15}$$

where $A_h := \pi_0^h$. Let $w := \sqrt{2\pi} R_g \pi$ (whether "$\pi$" means the circular ratio or stationary distribution should be clear), or

$$\pi(x_j) = \frac{1}{\sqrt{2\pi}} \exp(-(jh)^2) = \frac{1}{\sqrt{2\pi}} w_j \tag{3.16}$$

As $j \to \infty$, the leading behavior of $v_j$ is like

$$v_j = \exp\Big(-\sum_{k=1}^{j} \ln(1 + 2kh^2)\Big) \sim \exp(-C_h j \ln j)$$

which decays slower than $w_j$.

Clearly, $v_j$ is decreasing and

$$\sum_{k \geq j+1} h v_k \leq v_j h \sum_{m=1}^{\infty} \frac{1}{(1 + 2jh^2)^m} = \frac{v_j}{2jh}.$$

Hence, we find

$$\Big|\sum_{j \in \mathbb{Z}} h v_j - \sum_{j \in \mathbb{Z}} h w_j\Big| \leq \Big|\sum_{|j| \leq M} h|v_j - w_j|\Big| + \frac{2v_M}{2Mh} + \frac{C\pi(x_M)}{Mh}. \tag{3.17}$$

Moreover, since $-x \leq -\ln(1+x) \leq -x + \frac{1}{2}x^2$, using $\sum_{k=1}^{j} k^2 \leq j^3$, we find

$$w_j \exp(-jh^2) \leq v_j \leq w_j \exp(-jh^2 + 2j^3h^4).$$

It follows that $|v_j - w_j| \leq w_j C \max(jh^2, 2j^3h^4)$ for $j^3 h^4 \leq 1$ and $jh^2 \leq 1$. Since there exists independent of $h$ such that

$$\sum_{j \in \mathbb{Z}} h w_j \max(jh, 2j^3h^3) < C,$$

we find (3.17) can be controlled by

$$\Big|\sum_{j \in \mathbb{Z}} h v_j - \sum_{j \in \mathbb{Z}} h w_j\Big| \leq Ch + \Big(\frac{2v_M}{2Mh} + \frac{C\pi(x_M)}{Mh}\Big)|_{M=h^{-4/3}} \leq C_1 h.$$

Hence, $|h^{-1}A_h - \frac{1}{\sqrt{2\pi}}| \leq C_2 h$. Consequently, $h^{-1}\pi^h - R_g \pi$ is controlled by $h$ both in $\ell^\infty$ and in $\ell^1$.

## 3.2  Uniform error estimates

Note $\tilde{b}(x) = b(x)$ (since for $d = 1$ the detailed balance condition is satisfied always). We now use the equation for $q$ to investigate the uniform approximation of upwind scheme to the Fokker-Planck equation.

In [33, sect. 3.1], the following exponential decay has been proved:

**Proposition 3.1.** *Suppose that Assumption 2.1 holds and that the derivatives of $b$ and $\sigma$ are bounded. Then for any index $n > 0$, there exist a polynomial $p_n$ and $\gamma_n > 0$ such that*

$$\Big|\frac{\partial^n}{\partial x^n}(q(x,t) - 1)\Big| \leq p_n(x) \exp(-\gamma_n t).$$

Proposition 3.1, together with Lemma 2.1, implies that

**Theorem 3.2.** *Suppose that Assumption 2.1 holds and that the derivatives of $b$ and $\sigma$ are bounded. Then, for any $n \geq 0$, there exist $C_n > 0$ and $\tilde{\gamma}_n > 0$ such that*

$$\Big|\frac{\partial^n}{\partial x^n}(p(x,t) - \pi(x))\Big| \leq C_n \exp(-\nu_n|x|^2) \exp(-\tilde{\gamma}_n t). \tag{3.18}$$

*Suppose $\pi^h$ is the stationary solution for (3.11) with $\sum_j \pi_j^h = 1$ and recall $R_g$ (3.4). Then*

$$\sup_{t \geq 0} \sum_{j \in \mathbb{Z}} |p(x_j, t)h - p_j(t)| \leq Ch + 2\left\|R_g\pi - \frac{1}{h}\pi^h\right\|_{\ell^1} \leq Ch. \tag{3.19}$$

*Hence, for the upwind scheme (3.5), we have*

$$\sup_{t \geq 0} \|R_g\rho(\cdot, t) - \rho^h(t)\|_{\ell^1} \leq Ch. \tag{3.20}$$

*Proof.* Note that $p - \pi = \pi(q - 1)$. Lemma 2.1 and Proposition 3.1 imply (3.18).

We insert $\psi := p - \pi$ into the discrete Fokker-Planck equation (3.11) and by the standard Taylor expansion in numerical analysis scheme, we have

$$\frac{d}{dt}\psi(x_j, t) = \mathcal{L}_h^* \psi(x_j, t) + g(x_j, t)h,$$

where $\|g(x_j, t)\|_{\ell^1} \leq C \exp(-\gamma t)$ holds uniformly for small $h$ by (3.18). Then, we have $p(x_j, t) - \pi(x_j) = e^{t\mathcal{L}_h^*}(p(x_j, 0) - \pi(x_j)) + h \int_0^t e^{(t-s)\mathcal{L}_h^*} g \, ds$. Since the $p_j(t) = (e^{t\mathcal{L}_h^*} p^h(0))_j$ and $\pi^h = e^{t\mathcal{L}_h^*} \pi^h$, we have

$$p(x_j, t) - \pi(x_j) - \frac{1}{h}(p_j(t) - \pi_j^h) =$$
$$e^{t\mathcal{L}_h^*}\left(p(x_j, 0) - \pi(x_j) - \frac{1}{h}(p_j(0) - \pi_j^h)\right) + h \int_0^t e^{(t-s)\mathcal{L}_h^*} g(x_j, s) \, ds.$$

Since $e^{t\mathcal{L}_h^*}$ is $\ell^1$ contraction by Lemma 3.1, we have

$$\left\|R_g p(\cdot, t) - \frac{1}{h}p^h(t)\right\|_{\ell^1} \leq 2\left\|R_g \pi - \frac{1}{h}\pi^h\right\|_{\ell^1} + \left\|R_g p(\cdot, 0) - \frac{1}{h}p^h(0)\right\|_{\ell^1} + hC \int_0^t \exp(-\gamma z) dz.$$

The first claim (3.19) follows by noticing (3.8). The claim (3.20) follows from the relation between $\rho^h$ and $p^h$, using again (3.8). □

# 4 Properties of the jump process

We will investigate the forward and backward equations associated with the jump process $Z(t)$ corresponding to (3.11). It is convenient to introduce the Green's function

$$p_t(i, j) := \mathbb{P}(Z(t) = j | Z(0) = i) \geq 0. \tag{4.1}$$

Following [10, Chapter 2], we introduce the $Q$ matrix as

$$Q(i, j) = \frac{d}{dt} p_t(i, j)|_{t=0}. \tag{4.2}$$

## 4.1 Forward and backward equations

The Green's function $p_t(i, j)$ is a solution to (3.11) with the initial distribution $p_0(i, j) = \delta_{ij}$ (possibly not unique without Assumption 2.2). The equation for the Green's function is

$$\frac{d}{dt}p_t(i, j) = \alpha_{j-1} p_t(i, j-1) + \beta_{j+1} p_t(i, j+1) - (\alpha_j + \beta_j) p_t(i, j). \tag{4.3}$$

It follows that

$$Q(j, j) = -(\alpha_j + \beta_j), \ Q(j, j-1) = \beta_j, \ Q(j, j+1) = \alpha_j. \tag{4.4}$$

Recall the definition of irreducibility

**Definition 4.1.** *[10, Definition 2.47] A Markov chain is irreducible if $p_t(i, j) > 0$ for all $i, j$ and $t > 0$*

**Lemma 4.1.** *The jumping process $Z(t)$ corresponding to (3.11) is irreducible.*

*Proof.* First of all, we have

$$\frac{d}{dt}p_t(j, j) \geq -(\alpha_j + \beta_j)p_t(j, j),$$

which implies that $p_t(j, j) > 0$ for $t \geq 0$. Since

$$\frac{d}{dt}p_t(j, j-1) \geq \beta_j p_t(j, j) - (\alpha_{j-1} + \beta_{j-1})p_t(j, j-1),$$

$$\frac{d}{dt}p_t(j, j+1) \geq \alpha_j p_t(j, j) - (\alpha_{j+1} + \beta_{j+1})p_t(j, j+1),$$

it follows that $p_t(j, j-1) > 0$ and $p_t(j, j+1) > 0$ for $t > 0$.

In general, for $|k - j| =: m \geq 2$, we have by (4.3) that

$$p_t(j, k) = \int_0^t \exp(-(\alpha_j + \beta_j)(t-s))(\alpha_{k-1}p_s(j, k-1) + \beta_{k+1}p_s(j, k+1))\, ds.$$

By induction on $m$, we can prove that $p_t(j, k) > 0$ for $t > 0$. □

By [10, Corollary 2.58] and Lemma 3.2, if Assumption 2.2 holds, the chain is recurrent. The backward equation corresponding to the forward equation (3.11) reads

$$\frac{d}{dt}u_i(t) = \sum_{j\in\mathbb{Z}} Q(i, j)u_j(t) = \beta_i u_{i-1} - (\alpha_i + \beta_i)u_i + \alpha_i u_{i+1} =: (\mathcal{L}_h u)_i. \tag{4.5}$$

Clearly, $\mathcal{L}_h : \mathbb{R}^{\mathbb{Z}} \to \mathbb{R}^{\mathbb{Z}}$ is the dual operator of $\mathcal{L}_h^*$. In fact, letting

$$\langle u, v \rangle_h := \sum_{j\in\mathbb{Z}} h u_j v_j, \tag{4.6}$$

we have

$$\langle \mathcal{L}_h g, f \rangle_h = \langle g, \mathcal{L}_h^* f \rangle_h,$$

for any test sequence $f$ that has finite nonzero entries. (Note that sequences with finite nonzero entries are dense in $\ell^p$ with $p < \infty$, so this is general enough.) Let $u(t) = (u_j(t))_{j\in\mathbb{Z}}$ be the solution of (4.5). The semigroup defined by

$$e^{t\mathcal{L}_h}u(0) := u(t) \tag{4.7}$$

is the dual of $e^{t\mathcal{L}_h^*}$.

It is well-known that besides the forward equation (4.3), the Green's function also satisfies the backward equation (see [10, Theorem 2.14]):

$$\frac{d}{dt}p_t(i, j) = \sum_{k\in\mathbb{Z}} Q(i, k)p_t(k, j) = \beta_i p_t(i-1, j) - (\alpha_i + \beta_i)p_t(i, j) + \alpha_i p_t(i+1, j). \tag{4.8}$$

Formally, $P = e^{tQ}$ and we have $Qe^{tQ} = e^{tQ}Q$. This fact is an analogy to the continuous case (2.13). Since the chain is irreducible and recurrent, by [10, Corollary 2.34], the total probability is conserved $\sum_j p_t(i, j) = 1$ for all $i$ (i.e. no probability leaks to infinity). By [10, Theorem 2.26] and [10, Exercise 2.38], the backward equation (4.8) has a unique bounded solution in $\ell^\infty$ given any initial data $u(0) \in \ell^\infty$. Correspondingly, for general initial data $p^h(0) \in \ell^1$, the solution is a linear combination of $p_t(i, j)$. Hence, the forward equation is also well-posed, nonnegativity preserving and it preserves sum

$$\sum_{j\in\mathbb{Z}} p_j(t) = \sum_{j\in\mathbb{Z}} p_j(0). \tag{4.9}$$

Hence $e^{t\mathcal{L}_h}$ maps $\ell^\infty$ to $\ell^\infty$ and the semigroup $e^{t\mathcal{L}_h^*}$ given in (3.12) maps $\ell^1$ to $\ell^1$.

Note that though the Green's function $p_t(i, j)$ satisfies the backward equation, the probability distribution $p_i(t)$ for general initial data does not. Instead, the lemma below shows that $\sum_i p_t(j, i)u_i(0)$ satisfies the backward equation. Before we state the results, we introduce the weighted $\ell^p$ spaces here, which are analogies of the weighted $L^p(w)$ spaces in section 2.2. Given $w$ with $w_j \geq 0$, we define $\ell^p(w)$ as

$$\ell^p(w) := \left\{ q : \|q\|_{\ell^p(w)} := (\sum_{j\in\mathbb{Z}} w_j|q_j|^p)^{1/p} < \infty \right\}. \tag{4.10}$$

**Proposition 4.1.** *Let* $S(t) := e^{t\mathcal{L}_h}$. *Then,*

*(1) For any* $u(0) \in \ell^\infty$. *It holds that*

$$(S(t)u(0))_j = \sum_{i \in \mathbb{Z}} p_t(j, i) u_i(0). \tag{4.11}$$

*(2) The semigroup* $S(t)$ *is TVD, i.e., if* $u(0) \in \ell^1$, *then* $\sum_j |u_j(t) - u_{j-1}(t)|$ *is nonincreasing.*

*(3)* $S(t)$ *is symmetric in* $\ell^2(\pi^h)$ *for any* $t \geq 0$.

*(4)* $S(t)$ *is a contraction in* $\ell^p(\pi^h)$ *for any* $p \in [1, \infty]$.

*Proof.* (1). Let $v_j(t) = \sum_i p_t(j, i) u_i(0)$. Using Fubini theorem, we find that $v \in \ell^\infty$. Moreover, since $p_t(j, \cdot) \in \ell^1$ for all $j$ and $t \geq 0$, we find by (4.8) that

$$\frac{d}{dt} v_j(t) = \sum_{i \in \mathbb{Z}} (\beta_j p_t(j-1, i) + \alpha_j p_t(j+1, i) - (\alpha_j + \beta_j) p_t(j, i)) u_i(0)$$

$$= \beta_j v_{j-1}(t) + \alpha_j v_{j+1}(t) - (\alpha_j + \beta_j) v_j(t)$$

Hence, $v = u$ by the uniqueness of the bounded solution.

(2). The backward equation (4.5) can be rearranged into $\frac{d}{dt} u_j = \alpha_j(u_{j+1} - u_j) - \beta_j(u_j - u_{j-1})$. It follows that

$$\frac{d}{dt}(u_{j+1} - u_j) = \alpha_{j+1}(u_{j+2} - u_{j+1}) - (\alpha_j + \beta_{j+1})(u_{j+1} - u_j) + \beta_j(u_j - u_{j-1}).$$

This is a forward equation for the sequence $\{u_{j+1} - u_j\}$ and the rates are given so that the equation is well-posed. Note that $\{u_j(0) - u_{j-1}(0)\} \in \ell^1$ since $u(0) \in \ell^1$. Since well-posed forward equations are $\ell^1$ contractions, $S(t)$ is TVD. (Intuitively, we can multiply $\sigma_j := \text{sgn}(u_{j+1} - u_j)$ on both sides of the equations and use $\sigma_j(u_{j+2} - u_{j+1}) \leq |u_{j+2} - u_{j+1}|$, $\sigma_j(u_j - u_{j-1}) \leq |u_j - u_{j-1}|$ to obtain

$$\frac{d}{dt}|u_{j+1} - u_j| \leq \alpha_{j+1}|u_{j+2} - u_{j+1}| - (\alpha_j + \beta_{j+1})|u_{j+1} - u_j| + \beta_j|u_j - u_{j-1}|.)$$

(3). We denote $S := S(1)$ and $p(i, j) := p_1(i, j)$. Clearly, we only have to show that $S$ is symmetric by the semigroup property. Using the detailed balance, we have:

$$\sum_j \pi_j^h f_j (Sg)_j = \sum_j \sum_i f_j g_i \pi_j^h p(j, i) = \sum_{ij} \pi_i^h p(i, j) f_j g_i = \sum_i \pi_i^h g_i (Sf)_i.$$

(4). Let $(u_j^i(t))$, $i = 1, 2$ be two solutions and define $\tilde{u}_j = u_j^1 - u_j^2$. Then $(\tilde{u}_j)$ is also a solution and for any convex function $\varphi$ it holds that

$$\frac{d}{dt} \varphi(\tilde{u}_j) = \mathcal{L}_h \varphi(\tilde{u}_j)_j + \alpha_j(\varphi(\tilde{u}_j) + \varphi'(\tilde{u}_j)(\tilde{u}_{j+1} - \tilde{u}_j) - \varphi(\tilde{u}_{j+1}))$$

$$+ \beta_j(\varphi(\tilde{u}_j) + \varphi'(\tilde{u}_j)(\tilde{u}_{j-1} - \tilde{u}_j) - \varphi(\tilde{u}_{j-1})) \leq \mathcal{L}_h \varphi(\tilde{u}_j)_j. \tag{4.12}$$

If $\varphi$ is not differentiable at $\tilde{u}_j$, $\varphi'(\tilde{u}_j)$ is understood as one element in the subdifferential. Multiplying $\pi_j^h$ and applying the detailed balance (3.13), we have $\frac{d}{dt} \pi_j^h \varphi(\tilde{u}_j) \leq \mathcal{L}_h^\star(\pi \varphi(\tilde{u}))_j$. Taking sum on $j$ yields that $\frac{d}{dt} \sum_j \pi_j^h \varphi(\tilde{u}_j) \leq 0$. Choosing $\varphi(z) = |z|^p$ which is convex, we have the claims for $p \in [1, \infty)$.

For $p = \infty$, we multiply $\sigma_j := \text{sgn}(\tilde{u}_j)$ on both sides of the equation and obtain

$$\frac{d}{dt}|\tilde{u}_j| \leq \mathcal{L}_h|\tilde{u}|_j.$$

This implies that $\|\tilde{u}\|_{\ell^\infty}$ is non-increasing. $\qquad \square$

An important observation is that the discrete scheme always satisfies the detailed balance. If we define

$$q^h(t) := (q_j(t))_{j \in \mathbb{Z}}, \ q_j(t) = \frac{p_j(t)}{\pi_j^h}, \tag{4.13}$$

then $q^h$ satisfies the backward equation using the detailed balance condition (3.13):

$$\frac{d}{dt} q_j = \beta_j q_{j-1} + \alpha_j q_{j+1} - (\alpha_j + \beta_j) q_j. \tag{4.14}$$

With this interpretation, the relation (4.11) can be checked directly:

$$q_j(t) = \frac{1}{\pi_j^h} \sum_{i \in \mathbb{Z}} p_i(0) p_t(i, j) = \sum_{i \in \mathbb{Z}} p_t(i, j) \frac{p_i(0)}{\pi_j^h}.$$

Using the detailed balance (3.13), we have $\pi_i^h p_t(i, j) = p_t(j, i) \pi_j^h$. Hence,

$$q_j(t) = \sum_{i \in \mathbb{Z}} p_t(j, i) q_i(0).$$

## 4.2 Convergence for the weak confinement

The theory for irreducible time continuous Markov chain with countable state space is well-developed. See [10, Chapter 2]. We now use these theories to establish some basic properties of the jump processes and the upwind scheme. We have the following:

**Proposition 4.2.** *Suppose Assumption 2.2 holds. The jump process $Z(t)$ for (3.11) satisfies*

$$p_t(i, j) \to \pi_j^h, \ t \to \infty, \ \text{for all } i, j.$$

*Moreover, if we assume $p_j(0) = \frac{h \rho_0(x_j)}{\|\rho_0\|_{\ell^1}} \leq C \pi_j^h$ for all $j \in \mathbb{Z}$, we then have*

$$\sum_{j \in \mathbb{Z}} |p_j(t) - \pi_j^h| \to 0, \ t \to \infty.$$

*Consequently, the upwind scheme (3.5) satisfies $\left\| \rho^h(t) - \frac{1}{h} \pi^h \|\rho_0\|_{L^1(\mathbb{R})} \right\|_{\ell^1} \to 0.$*

*Proof.* By [10, Theorem 2.88,Theorem 2.66], we have for all $i, j$ that $p_t(i, j) \to \pi_j^h$ as $t \to \infty$. Now, in general, we have

$$p_j(t) = \sum_{i \in \mathbb{Z}} p_i(0) p_t(i, j).$$

Since $|p_t(i, j)| \leq 1$, the dominant convergence theorem implies that

$$p_j(t) \to \pi_j^h, \ t \to \infty, \ \forall j \in \mathbb{Z}.$$

Equation (4.14) has the maximal principle following the last claim in Proposition 4.1:

$$|q_j(t) - \theta| \leq \max_{j \in \mathbb{Z}} |q_j(0) - \theta|, \ \forall \theta \in \mathbb{R}.$$

In particular, we take $\theta = 1$. By the assumption, we have $|q_j(0)| \leq C$ and thus $|q_j(t) - 1| \leq C_1, \ \forall t \geq 0$. Since $p_j(t) \to \pi_j^h$, we have $q_j(t) \to 1, \forall j$. Dominant convergence theorem then yields

$$\sum_{j \in \mathbb{Z}} \pi_j^h |q_j(t) - 1| = \sum_{j \in \mathbb{Z}} |p_j(t) - \pi_j^h| \to 0, \ t \to \infty.$$

Using the relation between $p^h$ and $\rho^h$, we find

$$\left\| \rho^h(t) - \frac{1}{h} \pi^h \|\rho_0^h\|_{\ell^1} \right\|_{\ell^1} \to 0, \ t \to \infty.$$

By equation (3.8), the claim follows. □

The above proof makes use of the boundedness of $p_t(i, j)$ heavily. This clearly has no correspondence in the continuous case as $h \to 0$. Naturally, one may wonder whether we have the convergence uniform in $h \to 0$. We will investigate this in the next section.

# 5 Large time behaviors for strong confinement

In section 4.2, we have seen that the distribution of the jump process converges to the stationary solution under the weak confinement assumption. However, we do not have any rate for the convergence. Under the strong confinement (Assumption 2.1), we know that the convergence of the distribution for SDE (2.1) in $L^1(\mathbb{R})$ norm is exponential, which is obtained by using relative entropy and log Sobolev inequality [29]. Naturally, we desire that under Assumption 2.1 the jump process (3.11) has uniform geometric ergodicity under $\ell^1$ norm.

The convergence of $p^h(t)$ to $\pi^h$ in total variation norm (or $h^{-1}p^h(t) \to h^{-1}\pi^h$ in $\ell^1$) is equivalent to convergence of $q^h(t)$ to 1 in $\ell^1(\pi^h)$. Hence, we can consider the geometric convergence of $q^h(t)$ to 1 in $\ell^p(\pi^h)$ $(p \geq 1)$, which is closely related to spectral gaps of the semigroup $\{e^{t\mathcal{L}_h}\}$. This is a typical Perron-Frobenius type question. Besides the traditional compactness requirement of the semigroup $\{e^{t\mathcal{L}_h}\}$ in $\ell^p(\pi^h)$, some sufficient conditions for the Perron-Frobenius type theorems include the hypercontractivity and uniform integrability [34, 27, 35]. The classical result of Gross [36] tells us that the hypercontractivity is equivalent to log-Sobolev inequality. Proving such type of results for processes with infinite discrete states is often difficult. Our strategy is to prove a discrete Poincaré inequality, which uses a different Lyapunov function compared with the log-Sobolev inequality.

In section 5.1, we look at the general Lyapunov functions and derive the energy dissipation relations. Then, we use the quadratic function as the Lyapunov function and derive the discrete Poincaré inequality. In section 5.2, we establish the uniform geometric ergodicity.

## 5.1 A discrete Poincaré inequality

Slightly different from equation (4.12), we note the following for a smooth function $\varphi$:

$$\frac{d}{dt}\varphi(q_j) = \mathcal{L}_h(\varphi'(q)q)_j + \beta_j q_{j-1}(\varphi'(q_j) - \varphi'(q_{j-1})) + \alpha_j q_{j+1}(\varphi'(q_j) - \varphi'(q_{j+1})). \quad (5.1)$$

By the detailed balance condition (3.13), this gives for convex function $\varphi$ that

$$\frac{d}{dt}\sum_{j\in\mathbb{Z}}\pi_j^h\varphi(q_j) = -\sum_{j\in\mathbb{Z}}\alpha_j\pi_j^h(q_j - q_{j+1})(\varphi'(q_j) - \varphi'(q_{j+1})) \leq 0. \quad (5.2)$$

This is the energy dissipation relation. If $\varphi(q) = q\log q - q + 1$, $\sum_j \pi_j^h\varphi(q_j)$ gives the relative entropy. What we find useful is the quadratic function $\varphi(q) = \frac{1}{2}(q - \sum_k \pi_k^h q_k)^2$. Then, we have

$$\frac{d}{dt}\mathcal{F}_h = -\mathcal{D}_h \quad (5.3)$$

with

$$\mathcal{F}_h := \frac{1}{2}\sum_{j\in\mathbb{Z}}\pi_j^h(q - \sum_{k\in\mathbb{Z}}\pi_k^h q_k)^2, \quad \mathcal{D}_h := \sum_{j\in\mathbb{Z}}\alpha_j\pi_j^h(q_j - q_{j+1})^2. \quad (5.4)$$

Now we need to control $\mathcal{F}_h$ using $\mathcal{D}_h$. This type of control is achieved by Poincaré inequality. Below is a lemma modified from [11, Proposition 1] or [35, Lemma 1.3.10], which is a discrete Hardy inequality. For the convenience of the readers, we also attach the proof in Appendix A.

**Lemma 5.1.** *Let $\theta$ be a non-negative sequence with $\sum_j \theta_j < \infty$ and $\mu$ be a positive sequence on $\mathbb{Z}$. Set*

$$A := \sup_f \left\{ \max\left(\sum_{j\geq 0}\theta_j(\sum_{k=0}^{j}f_k)^2, \sum_{j\leq -1}\theta_j(\sum_{k=j}^{-1}f_k)^2\right) : \sum_{j\in\mathbb{Z}}\mu_j f_j^2 = 1 \right\} \quad (5.5)$$

*and*

$$B := \max\left(\sup_{j\geq 0}(\sum_{k=0}^{j}\mu_k^{-1})\sum_{k\geq j}\theta_k \ , \ \sup_{j<0}(\sum_{k=j}^{-1}\mu_k^{-1})\sum_{k\leq j}\theta_k\right). \tag{5.6}$$

*Then it holds that $B \leq A \leq 4B$.*

Using Lemma 5.1 and the approach in [35, sect. 1.3.3], it is straightforward to find:

**Lemma 5.2.** *Let $\alpha$ and $\beta$ be the rates in (3.6) for the jump process $Z(t)$. Define*

$$\kappa := \inf_f\left\{\sum_{j\in\mathbb{Z}}\alpha_j\pi_j^h(f_{j+1}-f_j)^2 : \sum_{j\in\mathbb{Z}}\pi_j^h f_j^2 = 1, \sum_{j\in\mathbb{Z}}\pi_j^h f_j = 0\right\}.$$

*Then we have*

$$\kappa^{-1} \leq 8\max\left(\sup_{j\geq 0}(\sum_{k=0}^{j}(\alpha_k\pi_k^h)^{-1})\sum_{k\geq j+1}\pi_k^h, \ \sup_{j\leq 0}(\sum_{k=j}^{0}(\beta_k\pi_k^h)^{-1})\sum_{k\leq j-1}\pi_k^h\right).$$

*Proof.* Consider $\theta$, $\mu$, $A$ and $B$ in Lemma 5.1. Let

$$A_1 := \sup_g\left\{\sum_{j\geq 0}\theta_j(\sum_{k=0}^{j}g_k)^2 + \sum_{j\leq -1}\theta_j(\sum_{k=j}^{-1}g_k)^2 : \sum_{j\in\mathbb{Z}}\mu_j g_j^2 = 1\right\}.$$

Then we have $A \leq A_1 \leq 2A$.

Clearly, for any sequence $g$ we can define a sequence $f$ such that

$$f_0 = 0, \ g_k = f_{k+1} - f_k$$

and this is a one-to-one correspondence. Then, we can rewrite $A_1$ in terms of $f$ as

$$A_1 = \sup_f\left\{\sum_{j\geq 0}\theta_j f_{j+1}^2 + \sum_{j\leq -1}\theta_j f_j^2 : \sum_{j\in\mathbb{Z}}\mu_j(f_{j+1}-f_j)^2 = 1, \ f_0 = 0\right\}. \tag{5.7}$$

It is clear that

$$A_1 = \sup_f\left\{\frac{\sum_{j\geq 0}\theta_j(f_{j+1}-f_0)^2 + \sum_{j\leq -1}\theta_j(f_j-f_0)^2}{\sum_{j\in\mathbb{Z}}\mu_j(f_{j+1}-f_j)^2} : \right.$$
$$\left. f\not\equiv const, \ \sum_{j\in\mathbb{Z}}\mu_j(f_{j+1}-f_j)^2 < \infty\right\}. \tag{5.8}$$

Now we define $\theta_j = \pi_{j+1}^h$ for $j \geq 0$ and $\theta_j = \pi_j^h$ for $j \leq -1$, and let $\mu_j = \alpha_j\pi_j^h$. Then, $A_1$ under this particular choice of $\theta$ and $\mu$ is

$$A_1 = \sup_f\left\{\frac{\sum_{j\in\mathbb{Z}}\pi_j^h(f_j-f_0)^2}{\sum_{j\in\mathbb{Z}}\alpha_j\pi_j^h(f_{j+1}-f_j)^2} : f\not\equiv const, \ \sum_{j\in\mathbb{Z}}\alpha_j\pi_j^h(f_{j+1}-f_j)^2 < \infty\right\}. \tag{5.9}$$

It is then straightforward to find

$$A_1^{-1} = \inf_f\left\{\sum_{j\in\mathbb{Z}}\alpha_j\pi_j^h(f_{j+1}-f_j)^2 : \sum_{j\in\mathbb{Z}}\pi_j^h(f_j-f_0)^2 = 1\right\}. \tag{5.10}$$

In fact, if all sequences with $\sum_{j\in\mathbb{Z}}\alpha_j\pi_j^h(f_{j+1}-f_j)^2 < \infty$, $f\not\equiv const$ satisfy $\sum_{j\in\mathbb{Z}}\pi_j^h(f_j-f_0)^2 < \infty$, then (5.10) is clear. If there exists $f$ such that $\sum_{j\in\mathbb{Z}}\alpha_j\pi_j^h(f_{j+1}-f_j)^2 < \infty$ but $\sum_{j\in\mathbb{Z}}\pi_j^h(f_j-f_0)^2 = \infty$, then $A_1 = \infty$. If this case happens, we can then take $\bar{f}^N =$ $A_N(f_i 1_{|i| \leq N})_{i \in \mathbb{Z}}$ with $A_N$ picked so that $\sum_j \pi_j^h (\tilde{f}_j^N - \tilde{f}_0^N)^2 = 1$. Then, $A_N \to 0$ and the infimum in (5.10) over $\tilde{f}_N$ is zero. Hence, (5.10) holds.

Using (5.10), we have

$$\kappa = \inf_f \left\{ \sum_{j \in \mathbb{Z}} \alpha_j \pi_j^h (f_{j+1} - f_j)^2 : \sum_{j \in \mathbb{Z}} \pi_j^h (f_j - \sum_k f_k \pi_k^h)^2 = 1 \right\} \geq A_1^{-1}.$$

This is because for $f \in \ell^2(\pi^h)$, the constant $c$ that minimizes $\inf_c \sum_{j \in \ell^2(\pi^h)} \pi_j^h (f_j - c)^2$ is the mean $c = \sum_k f_k \pi_k^h$. Hence, we conclude by Lemma 5.1 that

$$\kappa \geq \frac{1}{2} A^{-1} \geq \frac{1}{8} B^{-1}.$$

Using the detailed balance $\alpha_k \pi_k^h = \beta_{k+1} \pi_{k+1}^h$ for $k \leq -1$, we have

$$B = \max \left( \sup_{j \geq 0} (\sum_{k=0}^j (\alpha_k \pi_k^h)^{-1}) \sum_{k \geq j+1} \pi_k^h, \ \sup_{j \leq 0} (\sum_{k=j}^0 (\beta_k \pi_k^h)^{-1}) \sum_{k \leq j-1} \pi_k^h \right).$$

The claim then follows. □

**Lemma 5.3.** *Suppose $S_1 \leq \sigma^2 \leq S_2$ for $S_2 > S_1 > 0$ and $b$ is a smooth function. Then, fixing $R > 0$, we can find $C(R) > 0$ and $h_0 > 0$ such that*

$$\max_{0 \leq j \leq [R/h]+1} \pi_j^h \leq C(R) \min_{0 \leq j \leq [R/h]+1} \pi_j^h, \ \forall h \leq h_0.$$

*and that*

$$\max_{-[R/h]-1 \leq j \leq 0} \pi_j^h \leq C(R) \min_{-[R/h]-1 \leq j \leq 0} \pi_j^h, \ \forall h \leq h_0.$$

*Proof.* We only prove the claim for $0 \leq j \leq [R/h]+1$. The other case is similar. By the relation, we have

$$\pi_j^h = \pi_0^h \prod_{k=1}^j \frac{\alpha_{k-1}}{\beta_k} = \pi_0^h \prod_{k=1}^j \frac{s_{k-1}^+/h + \sigma_{k-1/2}^2/(2h^2)}{s_k^-/h + \sigma_{k-1/2}^2/(2h^2)}. \tag{5.11}$$

Hence, for $h$ small enough, we have

$$\pi_0^h \prod_{k=1}^j \frac{1}{1 + 2h|s(x_k)|/S_1} \leq \pi_j^h \leq \pi_0^h \prod_{k=1}^j \left( 1 + 2h \frac{|s(x_{k-1})|}{S_1} \right). \tag{5.12}$$

Using (5.12), we find

$$\frac{\max_{0 \leq j \leq [R/h]+1} \pi_j^h}{\min_{0 \leq j \leq [R/h]+1} \pi_j^h} \leq \prod_{k=1}^{[R/h]+1} \left( 1 + 2h \frac{|s(x_{k-1})|}{S_1} \right) \prod_{k=1}^{[R/h]+1} \left( 1 + 2h \frac{|s(x_k)|}{S_1} \right).$$

Note that $\prod_{k=1}^{[R/h]+1} \left( 1 + 2h \frac{|s(x_k)|}{S_1} \right) \leq \exp(\frac{2}{S_1} \sum_{k=1}^{[R/h]+1} h|s(x_k)|)$. The inside of the right hand side is the Riemann sum for the integral $\frac{2}{S_1} \int_0^{R+h} |s(x)| \, dx$. Hence, the right hand side is bounded by a number depending on $R$ when $h$ is small enough. Similarly, $\prod_{k=1}^{[R/h]+1} (1 + 2h \frac{|s(x_{k-1})|}{S_1}) \leq C_1(R)$. The claim follows. □

Now, we are able to conclude the discrete Poincaré inequality:

**Theorem 5.1.** *Suppose Assumption 2.1 holds with $S_1 \leq \sigma^2 \leq S_2$. Let $\pi^h$ be the stationary distribution of the jump process $Z(t)$. Then the discrete Poincaré inequality holds for measure $\pi^h$ when $h$ is small enough. In other words, there exist $h_0 > 0$ and $\kappa_1 > 0$ independent $h$ so that for any $f \in \ell^2(\pi^h)$, we have*

$$\kappa_1 \left( \sum_{j \in \mathbb{Z}} \pi_j^h f_j^2 - (\sum_{k \in \mathbb{Z}} \pi_k^h f_k)^2 \right) \leq \sum_{j \in \mathbb{Z}} \alpha_j \pi_j^h (f_{j+1} - f_j)^2, \tag{5.13}$$

*where $\alpha_j$ is the rate in (3.6).*

*Proof.* By (3.6), $\alpha_k^{-1} \leq 2S_1^{-1}h^2$ and $\beta_k^{-1} \leq 2S_1^{-1}h^2$, we have the following bound:

$$B_1 := \max\Big(\sup_{j \geq 0}(\sum_{k=0}^{j}(\alpha_k \pi_k^h)^{-1}) \sum_{k \geq j+1} \pi_k^h, \ \sup_{j \leq 0}(\sum_{k=j}^{0}(\beta_k \pi_k^h)^{-1}) \sum_{k \leq j-1} \pi_k^h\Big)$$
$$\leq 2S_1^{-1} \max(I_+, \ I_-),$$

where the numbers are given by

$$I_+ := \sup_{j \geq 0} h^2(j+1)(\max_{0 \leq k \leq j}(\pi_k^h)^{-1}) \sum_{k \geq j+1} \pi_k^h, \ I_- := \sup_{j < 0} h^2(|j|+1) \max_{j+1 \leq k \leq 0}(\pi_k^h)^{-1} \sum_{k \leq j} \pi_k^h.$$

Moreover, we can find $R > 0$ such that $s(x) = b(x) - \sigma(x)\sigma'(x) < -r|x|$ for $x > R$ and $s(x) = b(x) - \sigma(x)\sigma'(x) > r|x|$ for $x < -R$. Hence,

$$\alpha_k = \frac{\sigma_{k+1/2}^2}{2h^2}, \ \text{for } kh > R; \ \beta_k = \frac{\sigma_{k-1/2}^2}{2h^2}, \ \text{for } kh < -R. \tag{5.14}$$

By (5.14), we have for $j \geq [R/h] + 1 =: j^*$ that

$$\pi_{j+n}^h = \pi_j^h \prod_{i=1}^{n} \frac{\alpha_{j+i-1}}{\beta_{j+i}} = \pi_j^h \prod_{i=1}^{n} \frac{\sigma_{i+j-1/2}^2}{\sigma_{i+j-1/2}^2 + 2hs_{i+j}^-} \leq \pi_j^h \prod_{i=1}^{n} \frac{1}{1 + 2hs_{i+j}^-/S_2}, \ n \geq 1.$$

Hence, we have

$$\sum_{k \geq j+1} \pi_k^h \leq \pi_j^h \sum_{k \geq j+1} \frac{1}{(1 + 2rh^2(j+1)/S_2)^{k-j}} = \frac{S_2}{2r} \frac{\pi_j^h}{(j+1)h^2}, \tag{5.15}$$

where we have used $s_{i+j}^- \geq r(j+1)h$ for $i \geq 1$.

Let $K := \frac{S_2}{2r}$. If $0 \leq j \leq [R/h] = j^* - 1$, we have by (5.15) that

$$\sum_{k \geq j+1} \pi_k^h \leq (j^* - j) \max_{0 \leq k \leq j^*} \pi_k^h + K \frac{\pi_{j^*}^h}{(j^* + 1)h^2}.$$

Consequently, by Lemma 5.3,

$$h^2(j+1)(\max_{0 \leq k \leq j}(\pi_k^h)^{-1}) \sum_{k \geq j+1} \pi_k^h \leq ((R+h)^2 + K)C(R),$$

and the right hand side is uniformly bounded for $h \leq h_0$.

If $j \geq j^*$, using (5.15) again, we have

$$h^2(j+1)(\max_{0 \leq k \leq j}(\pi_k^h)^{-1}) \sum_{k \geq j+1} \pi_k^h \leq K(\min_{0 \leq k \leq j} \pi_k^h)^{-1} \pi_j^h \leq KC(R).$$

The last inequality holds because

$$\min_{0 \leq k \leq j} \pi_k^h = \min\Big(\min_{0 \leq k \leq j^*} \pi_k^h, \pi_j^h\Big).$$

Clearly, $\pi_j^h \leq \pi_{j^*}^h$. If $\pi_j^h \geq \min_{0 \leq k \leq j^*} \pi_k^h$, then $(\min_{0 \leq k \leq j} \pi_k^h)^{-1} \pi_j^h \leq (\min_{0 \leq k \leq j^*} \pi_k^h)^{-1} \pi_{j^*}^h \leq C(R)$ by Lemma 5.3. Otherwise, $(\min_{0 \leq k \leq j} \pi_k^h)^{-1} \pi_j^h = 1$. Hence, $I_+$ is bounded.

Similarly, $I_-$ can be bounded. Overall, $B_1$ is bounded by a constant $M$ depending on $R, r, S_1, S_2$ and $h_0$. Then, by Lemma 5.2, we have

$$\kappa \geq \frac{1}{8B_1} \geq \frac{1}{8M}.$$

Taking $\kappa_1 = 1/(8M)$ finishes the proof. $\qquad\square$

## 5.2 Uniform ergodicity

Recall that $\ell^1$ and $\ell^p(w)$ are defined in equation (1.5) and equation (4.10) respectively. Using Theorem 5.1, we are able to conclude that

**Theorem 5.2.** *Suppose Assumption 2.1 holds with $S_1 \leq \sigma^2 \leq S_2$. Consider the jump process $Z(t)$ corresponding to (3.11) and $q$ defined by (4.13). Then we have*

$$\left\| q^h(t) - \sum_{j\in\mathbb{Z}} \pi_j^h q_j \right\|_{\ell^2(\pi^h)} = \|q^h(t) - 1\|_{\ell^2(\pi^h)} \leq \|q^h(0) - 1\|_{\ell^2(\pi^h)} e^{-\kappa_1 t}.$$

*Consequently, $p^h(t)$ converges to $\pi^h$ exponentially fast in the total variation norm:*

$$\sum_{j\in\mathbb{Z}} |p_j(t) - \pi_j^h| \leq C \exp(-\kappa_1 t), \ \forall t > 0.$$

*The upwind scheme (3.5) satisfies $\|\rho^h(t) - \frac{1}{h}\|\rho_0^h\|_{\ell^1}\pi_j^h\|_{\ell^1} \leq C \exp(-\kappa_1 t)$.*

*Proof.* Recall the definition of $\mathcal{F}_h$ and $\mathcal{D}_h$ in (5.4). Then, by Theorem 5.1, we have

$$\frac{d}{dt}\mathcal{F}_h = -\mathcal{D}_h \leq -2\kappa_1 \mathcal{F}_h.$$

Noticing $\sum_j \pi_j^h q_j = \sum_j p_j = 1$ and $\mathcal{F}_h = \|q - \sum_j \pi_j^h q_j\|_{\ell^2(\pi^h)}^2$, the first claim follows.

By Hölder's inequality, it holds that

$$\sum_{j\in\mathbb{Z}} |p_j(t) - \pi_j^h| = \|q^h(t) - 1\|_{\ell^1(\pi^h)} \leq \|q^h(t) - 1\|_{\ell^2(\pi^h)} \leq C \exp(-\kappa_1 t).$$

Since

$$\rho_j(t) = \frac{1}{h}\|\rho_0^h\|_{\ell^1} p_j(t),$$

we then have

$$\left\| \rho^h(t) - \frac{1}{h}\|\rho_0^h\|_{\ell^1}\pi^h \right\|_{\ell^1} \leq \|\rho_0^h\|_{\ell^1} \sum_j |p_j(t) - \pi_j^h| \leq C \exp(-\kappa_1 t).$$

$\square$

Using the second claim of Theorem (5.2), we conclude the following property of the semigroup $e^{t\mathcal{L}_h^*}$:

**Corollary 5.1.** *Suppose that $v \in \ell^1$ and $\sum_j hv_j = 0$. Then,*

$$\left\| e^{t\mathcal{L}_h^*} v \right\|_{\ell^1} \leq C \exp(-\kappa_1 t).$$

*Proof.* Let $v^+ = \{v_j \vee 0\}$ and $v^- = \{-v_j \wedge 0\}$ so that $v = v^+ - v^-$. Let

$$p^1(t) := e^{t\mathcal{L}_h^*} \frac{hv^+}{\|v^+\|_{\ell^1}}, \ p^2(t) := e^{t\mathcal{L}_h^*} \frac{hv^-}{\|v^-\|_{\ell^1}}.$$

By Theorem (5.2), we have

$$\sum_{j\in\mathbb{Z}} |p_j^i(t) - \pi_j^h| \leq C_i \exp(-\kappa_1 t), \ i = 1, 2,$$

for some constants $C_i$.

Note that $\sum_j hv_j = 0$ implies $\|v^+\|_{\ell^1} = \|v^-\|_{\ell^1} = \frac{1}{2}\|v\|_{\ell^1}$. We have

$$\|e^{t\mathcal{L}_h^*}v\|_{\ell^1} = \sum_{j\in\mathbb{Z}}\left| \|v^+\|_{\ell^1}p_j^1(t) - \|v^-\|_{\ell^1}p_j^2(t) \right| = \frac{1}{2}\|v\|_{\ell^1}\sum_{j\in\mathbb{Z}}|p_j^1(t) - p_2^j(t)| \leq C \exp(-\kappa_1 t).$$

$\square$

Corollary 5.1 tells us that $e^{t\mathcal{L}_h^*}$ has a spectral gap in $\ell^1$. For any $v \in \ell^1$, we define the projection onto the space spanned by $\pi^h$ as

$$Pv := \left( \sum_{j \in \mathbb{Z}} h v_j \right) \left( \frac{1}{h} \pi^h \right). \tag{5.16}$$

Clearly, $Pv$ is invariant under $e^{t\mathcal{L}_h^*}$. Corollary 5.1 implies that if $v$ has no component in the direction of $\pi^h$, then $e^{t\mathcal{L}_h^*}v$ decays exponentially fast.

Now, we are able to conclude Theorem 3.1, i.e. bounding the error for approximating $\pi(x_j)$ using $\pi_j^h$. Note that for $j \in \mathbb{Z}$

$$\mathcal{L}_h^* \left( \pi(x_j) - \frac{1}{h} \pi_j^h \right) = \mathcal{L}_h^*(\pi(x_j)) = \tau_j h, \tag{5.17}$$

where $|\tau_j| \leq C$ and $\sum_j h |\tau_j| \leq C$ by direct Taylor expansion and Lemma 2.1. Intuitively, $P(R_g \pi - \frac{1}{h} \pi^h) = O(h)$, and $\mathcal{L}_h^*$ has a spectral gap in $\ell^1$. Hence, we may possibly invert $\mathcal{L}_h^*$ and obtain

$$\left\| R_g \pi - \frac{1}{h} \pi^h \right\|_{\ell^1} \leq Ch.$$

This understanding is not quite a rigorous proof. Below, we provide a rigorous proof.

*Proof of Theorem 3.1.* We have the following identity for operators from $\ell^1$ to $\ell^1$:

$$I = e^{t\mathcal{L}_h^*} + \int_0^t e^{(t-s)\mathcal{L}_h^*} \mathcal{L}_h^* \, ds. \tag{5.18}$$

In fact, for any $v \in \ell^1$ not does not depend on time, we set $f = \mathcal{L}_h^* v$. Then, $\frac{d}{dt} v + \mathcal{L}_h^* v = f$ implies that $v(t) = e^{t\mathcal{L}_h^*} v(0) + \int_0^t \exp((t-s)\mathcal{L}_h^*) f(s) \, ds$. Since we have assumed $v(t) \equiv v$, the identity is proved.

Now, we act the identity on $E_j = \pi(x_j) - \frac{1}{h} \pi_j^h$. Using equation (5.17), we have

$$E = e^{t\mathcal{L}_h^*} E + h \int_0^t e^{(t-s)\mathcal{L}_h^*} \tau \, ds,$$

where $\|\tau\|_{\ell^1} \leq C$. Since $\tau$ is in the range of $\mathcal{L}_h^*$, we therefore have (recall (4.6))

$$\sum_{j \in \mathbb{Z}} h \tau_j = \langle 1, \tau_j \rangle_h = \langle \mathcal{L}_h 1, E \rangle_h = 0,$$

by approximating $E$ with sequences that have finite nonzero entries. Moreover, we define

$$\bar{\pi}_j = \frac{1}{h} \int_{x_j - h/2}^{x_j + h/2} \pi(y) \, dy,$$

and have $\|\bar{\pi} - R_g \pi\|_{\ell^1} \leq C_1 h$. Applying Corollary 5.1, we have

$$\|E\|_{\ell^1} \leq \|e^{t\mathcal{L}_h^*}(\bar{\pi} - R_g \pi)\|_{\ell^1} + \lim_{t \to \infty} \|e^{t\mathcal{L}_h^*}(\bar{\pi} - h^{-1}\pi^h)\|_{\ell^1} + h \int_0^\infty C e^{-(t-s)\kappa_1 t} \, ds.$$

The second term is zero by Corollary 5.1 and the result follows. $\qquad \square$

# 6 Finite domain with periodic boundary condition

If the domain is finite with periodic boundary condition or we consider the problems on torus with length $L$

$$\mathbb{T} = \mathbb{R}/(L\mathbb{Z}), \tag{6.1}$$

many of the proofs above can be significantly simplified. The Wiener process $W$ is the standard Wiener process in $\mathbb{R}$ wrapped into $\mathbb{T}$. Hence, the generator and the Kolmogorov equations are unchanged. For SDEs on torus, one may refer to [37, 38].

We will assume generally the following.

*Assumption* 6.1. Assume $b, \sigma$ are smooth functions on $\mathbb{T}$ and $\sigma^2 \geq S_1 > 0$.

By [37], section 2], Assumption 6.1 implies that the SDE has a unique stationary measure with smooth density. In fact for $d = 1$, we can verify this directly. Letting $v(x) = \pi(x)\sigma^2(x)$ and $b_1(x) = b(x)/\sigma^2(x)$, the equation $L^*\pi = 0$ implies that

$$v(x) = \exp\Big(-\int_0^x b_1(y)\,dy\Big)\Big(v(0) + C\int_0^x \exp\big(\int_0^z b_1(y)\,dy\big)\,dz\Big). \tag{6.2}$$

Using $v(L) = v(0)$, we find

$$v(0) + C\int_0^L \exp\Big(\int_0^z b_1(y)\,dy\Big)\,dz = v(0)\exp\Big(\int_0^L b_1(y)\,dy\Big) > 0,$$

which determines $C$ uniquely. Since $\int_0^x \exp(\int_0^z b_1(y)\,dy)\,dz \leq \int_0^L \exp(\int_0^z b_1(y)\,dy)\,dz$, $v(x) > 0$ for all $x \in [0, L]$. Hence, we can normalize so that $\int_0^L \pi(x)\,dx = 1$.

Note that for the Fokker-Planck equation on torus, the corresponding jump process may not be reversible (the stationary distribution does not have detailed balance). The function $q(x, t) = p(x, t)/\pi(x)$ satisfies (2.16) and the modified SDE is given by

$$dY = (\frac{1}{\pi}\partial_x(\sigma^2\pi) - b)\,dt + \sigma\,dW. \tag{6.3}$$

As before, $\pi$ is also the stationary solution to the modified SDE, and (2.22) still holds. With this observation, we have

**Lemma 6.1.** *Suppose Assumption 6.1 holds. Then, let $u(x, t) = \mathbb{E}_x\varphi(X)$ for the SDE (2.1) or $u(x, t) = \mathbb{E}_x\varphi(Y)$ for the modified SDE (6.3) where $\varphi \in C^\infty(\mathbb{T})$. Then for any integer $k > 0$ we have have for some $\lambda_k > 0$ that*

$$\|u - \langle \pi, \varphi \rangle\|_{C^k(\mathbb{T})} \leq C\exp(-\lambda_k t). \tag{6.4}$$

*Consequently, for any index $n$, we can find $\gamma_n > 0$ such that*

$$\sup_{x \in \mathbb{T}}\Big|\frac{\partial^n}{\partial x^n}(\rho(x, t) - \pi(x))\Big| \leq C_n\exp(-\gamma_n t). \tag{6.5}$$

The proof of Lemma 6.1 follows closely [33, section 6.1.2], and we put it in Appendix B for convenience. This fact is also used in [38, H3].

For the discretization, we pick a positive integer $N$ and define

$$h = \frac{L}{N}, \; x_j = jh, \; 0 \leq j \leq N - 1. \tag{6.6}$$

If $j$ falls out of $[0, N]$, we wrap it back into $[0, N]$ using periodicity. (For example, $j = N + 2$ will be understood as $j = 2$.)

**Lemma 6.2.** *Equation (3.11) has on $\mathbb{T}$ has a unique stationary solution up to multiplicative constants. Besides, the one with $\sum_j \pi_j^h = 1$ satisfies $\pi_j^h > 0$ for all $j$. Moreover, we have for any sequence $f$ that*

$$-\sum_{j=0}^{N-1} \pi_j^h f_j \mathcal{L}_h f_j = \sum_{j=0}^{N-1} \frac{\beta_{j+1}\pi_{j+1}^h + \alpha_j\pi_j^h}{2}(f_{j+1} - f_j)^2, \tag{6.7}$$

*where $\mathcal{L}_h$ is the generator of the jump process $Z(t)$ for (3.11) on $\mathbb{T}$.*

*Proof.* Note that the jump process $Z(t)$ is irreducible and aperiodic with finite states. The existence of a unique stationary distribution follows from the standard theory of Markov chains. See [10], for example. This stationary distribution (denoted as $\pi^h$) is clearly a positive solution of $\mathcal{L}_h^* f = 0$ with $\sum_j \pi_j^h = 1$. We fix this $\pi^h$ now, and show that all solutions are multiples of $\pi^h$.

Direct computation shows that for any $j = 0, \ldots, N-1$

$$f_j \mathcal{L}_h f_j = \frac{1}{2} (\mathcal{L}_h f^2)_j - \frac{\beta_j}{2} (f_{j-1}^2 - f_j)^2 - \frac{\alpha_j}{2} (f_j - f_{j+1})^2.$$

Multiplying $\pi_j^h$ and taking the sum on $j$ yield (6.7).

According to (6.7), we find that $\mathcal{L}_h f = 0$ only has constant solutions. This means that the right eigenspace of $\mathcal{L}_h$ corresponding to eigenvalue 0 is one dimensional. Hence, the left eigenspace of $\mathcal{L}_h$ for eigenvalue 0 is also one dimensional. This means that $\mathcal{L}_h^* f = 0$ has a unique solution up to multiplying constants □

The stationary solution has the following property:

**Lemma 6.3.** *There exists a constant $C$ independent of $h$ such that for sufficiently small $h$*

$$\max_{0 \le j \le N-1} \pi_j^h \le C \min_{0 \le j \le N-1} \pi_j^h. \tag{6.8}$$

*Proof.* We introduce the variable

$$z_j := \pi_j^h / \pi(x_j), \ j = 0, \ldots, N-1.$$

Since $\pi(\cdot)$ is bounded from below and above, we only need to investigate $z_j$. Consider the equation for $\pi_j^h$. We thus have

$$-\left( \frac{s_j^+ \pi(x_j) z_j - s_{j-1}^+ \pi(x_{j-1}) z_{j-1}}{h} - \frac{s_{j+1}^- \pi(x_{j+1}) z_{j+1} - s_j^- \pi(x_j) z_j}{h} \right)$$
$$+ \frac{1}{2h^2} (\sigma_{j+1/2}^2 \pi(x_{j+1}) z_{j+1} - (\sigma_{j+1/2}^2 + \sigma_{j-1/2}^2) \pi(x_j) z_j + \sigma_{j-1/2}^2 \pi(x_{j-1}) z_{j-1}) = 0. \tag{6.9}$$

Since $\pi(x)$ is a solution to $\mathcal{L}^* = 0$, there exists $h_0 > 0$ such that for all $h \le h_0$,

$$\mathcal{L}_h^* \pi(x_j) = \tau_j h, \ \forall 0 \le j \le N-1, \tag{6.10}$$

where $\|\tau_j\|_{\ell^\infty} \le C_1$ uniformly for $h \le h_0$. Subtracting (6.9) with $z_j \mathcal{L}_h^* \pi(x_j)$ and using (6.10), we have

$$T_h z_j := -\left( s_{j-1}^+ \pi(x_{j-1}) \frac{z_j - z_{j-1}}{h} - s_{j+1}^- \pi(x_{j+1}) \frac{z_{j+1} - z_j}{h} \right)$$
$$+ \frac{1}{2h^2} (\sigma_{j+1/2}^2 \pi(x_{j+1})(z_{j+1} - z_j) - \sigma_{j-1/2}^2 \pi(x_{j-1})(z_j - z_{j-1})) = -z_j \tau_j h. \tag{6.11}$$

Expanding $\pi(x_{j\pm1})$ in $\sigma_{j\pm1/2}^2 \pi(x_{j\pm1})$ terms around $x_{j\pm1/2}$, it is not hard to see $T_h$ is a first order consistent difference scheme for the modified backward operator

$$\tilde{\mathcal{L}} q = \frac{1}{2} \partial_x (\pi \sigma^2 \partial_x q) - \left( \frac{1}{2} \sigma^2 \partial_x \pi - s\pi \right) \partial_x q, \tag{6.12}$$

which is clearly the same as the one in (2.20).

The crucial observation is that both $T_h$ and $\tilde{\mathcal{L}}$ with Dirichlet boundary conditions have maximum principles. This allows us to prove the stability of $T_h$. Let us now investigate this in detail. Assume $z_j$ attains the maximum value at $j^*$. Without loss of generality, we can assume $j^* = 0$. Then, define for $j = 0, \ldots, N-1$ that

$$\zeta_j := \frac{z_j}{z_0} - 1.$$

We find then

$$T_h \zeta_j = -\frac{z_j}{\|z\|_{\ell^\infty}} \tau_j h, \ \text{for } j = 1, \ldots, N-1,$$

$$\zeta_0 = \zeta_N = 0.$$

Consider the equation

$$\tilde{\mathcal{L}}\,\phi(x) = 1, \ \phi(0) = \phi(L) = 0.$$

By the maximum principle, $\phi(x) < 0$ for $x \in (0, L)$. Since $T_h$ is a consistent scheme for $\tilde{\mathcal{L}}$, for sufficiently small $h$, we have

$$T_h \phi(x_j) \geq 1/2, \ j = 1, \ldots, N-1.$$

Letting $\xi_j := 2\|\tau\|_\infty \phi(x_j) h - \zeta_j$, we have for $j = 1, \ldots, N-1$,

$$T_h(\xi)_j \geq 0$$

with $\xi_0 = \xi_N = 0$. This means $\xi_j \leq 0$ by maximum principle and hence

$$\zeta_j \geq 2\|\tau\|_\infty \phi(x_j) h.$$

Similarly, replacing $\zeta$ with $-\zeta$, we have $\zeta_j \leq -2\|\tau\|_\infty \phi(x_j) h$. This means

$$\max_{0 \leq j \leq N-1} |\zeta_j| = \max_{0 \leq j \leq N-1} |\frac{z_j}{z_0} - 1| \leq 2\|\tau\|_\infty \|\phi\|_\infty h.$$

Hence, for all $j = 0, \ldots, N-1$,

$$\frac{z_j}{z_0} \geq 1 - 2\|\tau\|_\infty \|\phi\|_\infty h \geq \frac{1}{2},$$

when $h$ is sufficiently small. The claim (6.8) follows since $\pi$ is bounded from above and below by positive numbers. $\qquad\square$

Now, we prove the uniform consistency of the upwind schemes, which is an analogy of Theorem 3.1 and Theorem 3.2.

**Theorem 6.1.** *Consider the upwind scheme* (3.5) *and the jump process* $Z(t)$ *corresponding to* (3.11) *on* $\mathbb{T}$. *Suppose Assumption* 6.1 *holds. Then,*

(i) *The stationary distribution of* (3.11) *satisfies that*

$$\max_{0 \leq j \leq N-1} \left| \frac{1}{h}\pi^h - \pi(x_j) \right| \leq Ch.$$

(ii) *The following uniform error estimate holds for* (3.5). $\sup_{t \geq 0} \|R_g \rho(\cdot, t) - \rho^h(t)\|_{\ell^1} \leq Ch.$

The first claim is essentially proved in the proof of Lemma 6.3. There, we have seen that $|z_j/\|z\|_\infty - 1| \leq Ch$. Since $|\sum_j h\pi(x_j) - 1| \leq C_1 h$ and $\sum_j z_j \pi(x_j) = 1$, we then conclude that $|h^{-1}\|z\|_{\ell^\infty} - 1| \leq C_2 h$. The second claim can be proved in the same way as in the proof of Theorem 3.2.

We now move on to the convergence to equilibrium for the upwind scheme. Using Lemma 6.3 and that the torus is a bounded domain, the following version of discrete Poincaré inequality (analogy of Theorem 5.1) can be proved in a straightforward way (one can refer to [12, Proposition 4.6] for similar discussion).

**Lemma 6.4.** *Suppose Assumption* 6.1 *holds. Then there exists* $h_0 > 0$ *and* $\kappa_1 > 0$, *so that for any sequence* $f$, *we have*

$$\kappa_1 \sum_{j=0}^{N-1} \pi_j^h \left( f_j - \sum_{i=0}^{N-1} \pi_i^h f_i \right)^2 \leq \sum_{j=0}^{N-2} \frac{\beta_{j+1}\pi_{j+1}^h + \alpha_j \pi_j^h}{2}(f_{j+1} - f_j)^2. \tag{6.13}$$

*Proof.* Since $f_j - f_0 = \sum_{k=1}^{j}(f_k - f_{k-1})$, we have

$$\sum_{j=0}^{N-1} \pi_j^h (f_j - f_0)^2 \leq \sum_{j=1}^{N-1} \pi_j^h j \sum_{k=1}^{j}(f_k - f_{k-1})^2$$

$$= \sum_{k=1}^{N-1} \frac{\beta_k \pi_k^h + \alpha_{k-1}\pi_{k-1}^h}{2}(f_k - f_{k-1})^2 \sum_{j \geq k} \frac{2j\pi_j^h}{\beta_k \pi_k^h + \alpha_{k-1}\pi_{k-1}^h}.$$

The claim follows from the fact that when $h$ is sufficiently small

$$\sum_{k \leq j \leq N-1} \frac{2j\pi_j^h}{\beta_k \pi_k^h + \alpha_{k-1}\pi_{k-1}^h} \leq \frac{2N^2}{\min_{j,k}(\beta_k \pi_k^h/\pi_j^h + \alpha_{k-1}\pi_{k-1}^h/\pi_j^h)}$$

$$\leq \frac{2CN^2}{\min_k(\beta_{k+1} + \alpha_k)}$$

$$\leq \frac{2C}{S_1}N^2 h^2,$$

where we have applied Lemma 6.3 to obtain $\min_j \pi_k^h/\pi_j^h \geq \frac{1}{C}$ and $\min_j \pi_{k-1}^h/\pi_j^h \geq \frac{1}{C}$ for any $k$. Since $Nh = L$ and $\sum_j \pi_j^h(f_j - \sum_i \pi_i^h f_i)^2 \leq \sum_j \pi_j^h(f_j - f_0)^2$, the claim follows. $\quad\square$

The chain in general is not reversible. In fact, for the stationary solutions, we have

$$J_{j+1/2} = J = const.$$

If $J = 0$, then we must have $\prod_{j=0}^{N-1} \alpha_j = \prod_{j=0}^{N-1} \beta_j$, which may not be true. Hence, in general $J \neq 0$ and the process is not reversible. Defining

$$\tilde{\beta}_j := \frac{\alpha_{j-1}\pi_{j-1}^h}{\pi_j^h}, \ \tilde{\alpha}_j := \frac{\beta_{j+1}\pi_{j+1}^h}{\pi_j^h}, \ j = 0, \dots, N-1 \tag{6.14}$$

we have

$$\alpha_j + \beta_j = \tilde{\alpha}_j + \tilde{\beta}_j, \ j = 0, \dots, N-1.$$

Hence, using (3.11), we can write the equation for $q^h = p^h/\pi^h$ ($p^h$ and $q^h$ are similarly defined as in (3.10) and (4.13)) as

$$\frac{d}{dt}q_j = \tilde{\beta}_j q_{j-1} + \tilde{\alpha}_j q_{j+1} - (\tilde{\alpha}_j + \tilde{\beta}_j)q_j =: (\tilde{\mathcal{L}}_h q^h)_j, \ j = 0, \dots, N-1. \tag{6.15}$$

It is easily verified that $\pi^h$ is also a stationary solution of $\tilde{\mathcal{L}}_h^*$, the dual operator of $\tilde{\mathcal{L}}_h$:

$$(\tilde{\mathcal{L}}_h^* \pi^h)_j = \tilde{\alpha}_{j-1}\pi_{j-1}^h - (\tilde{\alpha}_j + \tilde{\beta}_j)\pi_j^h + \tilde{\beta}_{j+1}\pi_{j+1}^h = \beta_j \pi_j^h - (\alpha_j + \beta_j)\pi_j^h + \alpha_j \pi_j^h = 0. \tag{6.16}$$

With the preparation, we easily conclude the following, similar to Theorem 5.2.

**Theorem 6.2.** *Consider the upwind scheme* (3.5) *and the equivalent discrete Fokker-Planck equation* (3.11) *on torus. Suppose Assumption 6.1 holds. Then, we have* $\|q^h(t) - 1\|_{\ell^2(\pi^h)} \leq \|q^h(0) - 1\|_{\ell^2(\pi^h)}e^{-\kappa_1 t}$. *Consequently,* $p^h(t)$ *converges to* $\pi^h$ *exponentially fast in total variation norm* $\sum_j |p_j(t) - \pi_j^h| \leq C\exp(-\kappa_1 t)$, *and thus* $\|\rho^h(t) - \frac{1}{h}\|\rho_0^h\|_{\ell^1}\pi^h\|_{\ell^1} \leq C\exp(-\kappa_1 t)$.

*Proof.* Let $\varphi$ be a smooth function defined on $\mathbb{T}$. Applying (6.15) and using similar calculation as in equation (4.12), we have

$$\frac{d}{dt}\sum_{j=0}^{N-1}\pi_j^h \varphi(q_j) = \sum_{j=0}^{N-1}\pi_j^h \tilde{\alpha}_j(\varphi(q_j) + \varphi'(q_j)(q_{j+1} - q_j) - \varphi(q_{j+1}))$$

$$+ \sum_{j=0}^{N-1}\pi_j^h \tilde{\beta}_j(\varphi(q_j) + \varphi'(q_j)(q_{j-1} - q_j) - \varphi(q_{j-1})). \tag{6.17}$$

If we take $\varphi(q_j) = \frac{1}{2}(q_j - \sum_i \pi_i^h q_i)^2$, we then have by (6.17) and (6.14) that

$$\frac{1}{2}\frac{d}{dt}\sum_{j=0}^{N-1}\pi_j^h(q_j - \sum_i \pi_i^h q_i)^2 = -\sum_{j=0}^{N-1}\frac{\tilde{\alpha}_j \pi_j^h + \tilde{\beta}_{j+1}\pi_{j+1}^h}{2}(q_{j+1} - q_j)^2$$

$$= -\sum_{j=0}^{N-1}\frac{\beta_{j+1}\pi_{j+1}^h + \alpha_j \pi_j^h}{2}(q_{j+1} - q_j)^2 \tag{6.18}$$

$$\leq -\sum_{j=0}^{N-2}\frac{\beta_{j+1}\pi_{j+1}^h + \alpha_j \pi_j^h}{2}(q_{j+1} - q_j)^2.$$

Using Lemma 6.4, the remaining proof is similar to the proof of Theorem 5.2, and we omit. □

# 7 A Monte Carlo method for the upwind scheme

In this section, we propose some Monte Carlo methods to approximate the upwind scheme (3.5). One idea is to construct a jump process $\{Z_n^{\Delta t}\}$ with transition probability $\tilde{P} = I + \Delta t \, Q$ using forward Euler scheme in time. In other words, the probability distribution satisfies

$$p^{n+1} = (I + \Delta t \, Q)p^n, \tag{7.1}$$

where $p^n$ refers to the probability distribution at $n$-th step. There are two drawbacks. Firstly, the forward Euler introduces numerical errors in time discretization; secondly $I + \Delta t \, Q$ may have negative entries for any $\Delta t$. Another idea is to use the continuous time random walk. The process waits for a random time that satisfies an exponential distribution at a site and then performs a jump. This idea can avoid using the time discretization to recover (3.5). If we consider the upwind scheme on $\mathbb{R}$, we need the exponential distribution for the waiting time to depend on the site $j$, and a corresponding Monte Carlo method can be developed. For the jump process $Z(t)$ on torus, we can choose the exponential distribution independent of the sites. Then the number of jumps is a Poisson process and this motivates another Monte Carlo algorithm. For the convenience, we focus on the problems on torus only and explain this Monte Carlo algorithm in detail.

**Lemma 7.1** ([10, Example 2.5]). *Let $P$ be a transition matrix. Let $\mathcal{N}(t)$ be a Poisson process of intensity $\lambda$. If $Z_1(t)$ is the process that takes transitions at jumps of $\mathcal{N}(t)$ according to $P$, then $Z_1(t)$ is a continuous time jump process with $Q$ matrix to be*

$$Q = \lambda(P - I). \tag{7.2}$$

Recall that $Q$ matrix is defined in (4.2) so that $p_t(i,j) = \mathbb{P}(Z_1(t) = j | Z_1(0) = i)$ satisfies

$$\frac{d}{dt}p_t(i,j) = \sum_k Q(i,k)p_t(k,j) = \sum_k p_t(i,k)Q(k,j).$$

Lemma 7.1 follows easily from the fact $Z_1(t)$ is Markovian and that

$$p_t(i,j) = e^{-\lambda t} \sum_{n=0}^{\infty} \frac{(\lambda t)^n}{n!} P^n(i,j).$$

Here, $P^n$ is defined inductively by $P^{m+1}(i,j) = \sum_k P^m(i,k)P(k,j)$ with $P^1 = P$. With Lemma 7.1, we find that if $Q(i,j)$ is bounded, we can take $\lambda$ large enough so that

$$P = I + \lambda^{-1}Q \tag{7.3}$$

has nonnegative entries. Such construction is possible for problems on torus so that $Z_1(t)$ is a realization of $Z(t)$. This then gives the following Monte Carlo method for the upwind scheme:

1. Fix $T > 0$. Pick $\lambda \geq \max(\alpha + \beta)$ with $\alpha, \beta$ in (3.6). Pick $M$ for the number of samples.

2. For $m = 1 : M$:

   - Sample $\mathcal{N} \sim Poisson(\lambda T)$, and $j_0 \sim p_j(0)$.
   - Sample $Y_{\mathcal{N}}$ according to the $j_0$-th row of $P^{\mathcal{N}}$. (In other words, we have a discrete time Markov chain $\{Y_n\}_{n=1}^{\mathcal{N}}$ with $Y_0 = j_0$ and transition matrix $P$ in (7.3), or $P(j,j) = 1 - \lambda^{-1}(\alpha_j + \beta_j)$, $P(j,j-1) = \lambda^{-1}\beta_j$, and $P(j,j+1) = \lambda^{-1}\alpha_j$.)

3. Let $\tilde{p}$ be the empirical distribution of $Y_{\mathcal{N}}$ (with $M$ values of $Y_{\mathcal{N}}$). Then, $\tilde{\rho}(x_j, T) = h^{-1}\|\rho_0^h\|_{\ell^1}\tilde{p}_j$ is the numerical solution.

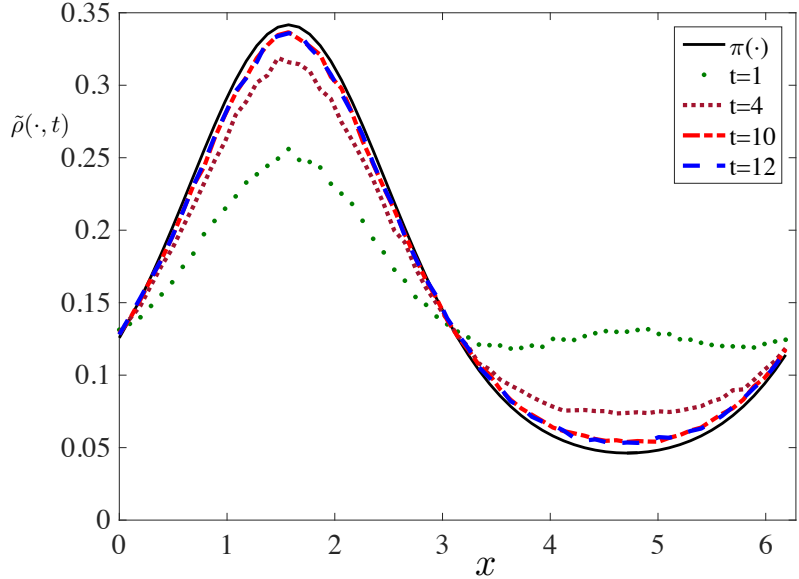

Figure 1: Monte Carlo simulation of the jump process. Number of grids $N = 2^6$, $\lambda \approx 291.7$ and number of samples $M = 10^6$. The solid black line shows the exact stationary solution $\pi(\cdot)$. Others show the computed numerical solution at $t = 1$(green dots), $t = 4$ (brown dotted line), $t = 10$ (red dash-dotted line) and $t = 12$ (blue dashed line). The stationary solution of the numerical solution is close to the stationary distribution of the SDE.

**Remark 7.1.** Since $\mathbb{E}\mathcal{N} = \lambda T$, $\lambda^{-1}$ is like the time step. Hence, $\lambda^{-1} \max(\alpha + \beta) \leq 1$ is like the CFL condition (for parabolic equations).

Note that we may use fast algorithms to pre-compute $P^n$ to save time. Consider the following SDE on $\mathbb{T}$ with $L = 2\pi$ and

$$b(x) = \cos(x)\exp(\sin(x)), \ \ \sigma(x) = \exp\left(\frac{1}{2}\sin(x)\right).$$

The function

$$s(x) = b(x) - \sigma(x)\sigma'(x) = \frac{1}{2}\cos x \exp(\sin x).$$

It can be check easily that the stationary solution is

$$\pi(x) \propto \exp(\sin(x)).$$

By the symbol "$\pi$" in this example, whether we mean the circular ratio or the stationary solution should be clear in the context.

Now, we take $\rho(x, 0) = \frac{1}{2\pi}$ so that $\lim_{t\to\infty} \rho(x, t) = \pi(x)$. The initial distribution for $j_0$ is therefore the uniform distribution. Figure 1 shows the computed $\tilde{\rho}$ at $t = 1, 4, 10, 12$, where we take number of grid points $N = 2^6$, $h = 2\pi/N$, $\lambda = \max(\alpha_j + \beta) + 10 \approx 291.7$ and the number of samples $M = 10^6$. We find that numerical solution of the Monte Carlo method for the jump process indeed converges to a stationary solution fast. Moreover, the stationary solution of the numerical solution is close to the stationary distribution of the SDE. This example therefore verifies our theory and the Monte Carlo method.

# Acknowledgements

The authors would like to thank Prof. Fengyu Wang for discussing uniform integrability and spectral gaps at Beijng Normal University. The work of J.-G. Liu is partially supported by KI-Net NSF RNMS11-07444, NSF DMS-1514826 and NSF DMS-1812573.

# A   Proof of Lemma 5.1

*Proof of Lemma* 5.1. Recall that $\theta$ is a non-negative sequence with $\sum_j \theta_j < \infty$ and $\mu$ is a positive sequence on $\mathbb{Z}$. We first pick $f_i = \mu_i^{-1} 1_{[0,M]}(i)$. By the definition of $A$, we have

$$A \sum_{k=0}^M \mu_k^{-1} = A \sum_{k=-\infty}^\infty \mu_k f_k^2 \geq \sum_{j \geq 0} \theta_j (\sum_{k=0}^j f_k)^2 \geq \sum_{j \geq M} (\sum_{k=0}^M \mu_k^{-1})^2 \theta_j.$$

Similarly, if we pick $f_i = \mu_i^{-1} 1_{[-M,-1]}(i)$, we have

$$A \sum_{k=-M}^{-1} \mu_k^{-1} = A \sum_{k=-\infty}^\infty \mu_k f_k^2 \geq \sum_{j \leq -1} \theta_j (\sum_{k=j}^{-1} f_k)^2 \geq \sum_{j \leq -M} (\sum_{k=-M}^{-1} \mu_k^{-1})^2 \theta_j.$$

This verifies that $A \geq B$.

On the other hand, let us assume $\sum_j \mu_j f_j^2 = 1$. Note the basic inequality

$$\frac{b-a}{2\sqrt{b}} \leq \sqrt{b} - \sqrt{a}, \ a \geq 0, b > 0. \tag{A.1}$$

Now let $\gamma_j := \sum_{k=0}^j \mu_k^{-1}$. Applying (A.1) and noting $\gamma_0 = \mu_0^{-1}$, we obtain

$$\sum_{k=0}^j \frac{\mu_k^{-1}}{\sqrt{\gamma_k}} = \frac{\mu_0^{-1}}{\sqrt{\gamma_0}} + \sum_{k=1}^j \frac{\gamma_k - \gamma_{k-1}}{\sqrt{\gamma_k}} \leq \sqrt{\gamma_0} + 2\sqrt{\gamma_j} - 2\sqrt{\gamma_0} \leq 2\sqrt{\gamma_j}. \tag{A.2}$$

Similarly,

$$\sum_{j \geq k} \frac{\theta_j}{\sqrt{\sum_{i \geq j} \theta_i}} = \sum_{j \geq k} \frac{\sum_{i \geq j} \theta_i - \sum_{i \geq j+1} \theta_i}{\sqrt{\sum_{i \geq j} \theta_i}} \leq 2\sqrt{\sum_{i \geq k} \theta_i}. \tag{A.3}$$

Consequently, we find

$$\begin{aligned}
\sum_{j \geq 0} \theta_j (\sum_{k=0}^j f_k)^2 &\leq \sum_{j \geq 0} \theta_j \Big( \sum_{k=0}^j f_k^2 \mu_k \sqrt{\gamma_k} \Big) \Big( \sum_{k=0}^j \frac{\mu_k^{-1}}{\sqrt{\gamma_k}} \Big) \\
&\leq 2 \sum_{j \geq 0} \theta_j \sqrt{\gamma_j} \sum_{k=0}^j f_k^2 \mu_k \sqrt{\gamma_k} \\
&\leq 2\sqrt{B} \sum_{j \geq 0} \frac{\theta_j}{\sqrt{\sum_{i \geq j} \theta_i}} \sum_{k=0}^j f_k^2 \mu_k \sqrt{\gamma_k} \\
&= 2\sqrt{B} \sum_{k \geq 0} f_k^2 \mu_k \sqrt{\gamma_k} \sum_{j \geq k} \frac{\theta_j}{\sqrt{\sum_{i \geq j} \theta_i}} \\
&\leq 4B \sum_{k \geq 0} f_k^2 \mu_k \leq 4B.
\end{aligned}$$

The first inequality is due to Hölder inequality. The second inequality is due to (A.2). The third inequality is due to (recall the definition of $\gamma_j$ and definition of $B$)

$$\sqrt{\gamma_j} \sqrt{\sum_{i \geq j} \theta_i} \leq \sqrt{B}.$$

The second last inequality is due to (A.3)

$$\sqrt{\gamma_k} \sum_{j \geq k} \frac{\theta_j}{\sqrt{\sum_{i \geq j} \theta_i}} \leq 2\sqrt{\gamma_k} \sqrt{\sum_{i \geq k} \theta_i} \leq 2\sqrt{B}.$$

Similarly, defining $\gamma_j = \sum_{k=j}^{-1} \mu_k^{-1}$, one can control

$$\sum_{j \leq -1} \theta_j \Big(\sum_{k=j}^{-1} f_k\Big)^2 \leq 4B.$$

Hence, $A \leq 4B$. $\qquad\qquad\square$

## B   Proof of Lemma 6.1

*Proof of Lemma* 6.1. Recall the notation

$$\langle \pi, f \rangle = \int_{\mathbb{T}} f(x)\pi(x)\,dx.$$

Without loss of generality, we assume $\langle \pi, \varphi \rangle = 0$ and consider the equation of $u$ for SDE (2.1) (the proof for the modified SDE (6.3) is just the same):

$$\partial_t u = \mathcal{L}u = b\,\partial_x u + \frac{1}{2}\Lambda \partial_{xx} u. \tag{B.1}$$

We see $\langle \pi, u \rangle = 0$ for all $t > 0$. Multiplying $2u$, we have

$$\partial_t |u|^2 = \mathcal{L}|u|^2 - \Lambda|\partial_x u|^2.$$

Multiplying $\pi$ and integrating yields

$$\frac{d}{dt}\int_{\mathbb{T}} \pi(x)|u|^2(x)\,dx = -\int_{\mathbb{T}} \pi\Lambda|\partial_x u|^2\,dx \leq -\lambda \int_{\mathbb{T}} \pi|u|^2\,dx. \tag{B.2}$$

The inequality follows from Poincaré inequality since $\langle \pi, u \rangle = 0$. We then obtain the exponential decay of $\langle \pi, |u|^2 \rangle$:

$$\int_{\mathbb{T}} |u|^2\,d\pi \leq \langle \pi, \varphi^2 \rangle \exp(-\lambda t).$$

Consequently, multiplying $e^{(\lambda-\delta)t}$ in (B.2) for $\delta > 0$ small and taking integral,

$$\frac{1}{2}\int_0^\infty e^{(\lambda-\delta)t}\int_{\mathbb{T}} \pi\Lambda|\partial_x u|^2\,dx = -\int_0^\infty e^{(\lambda-\delta)t}\frac{d}{dt}\int_{\mathbb{T}} \pi|u|^2\,dx dt \leq C.$$

This means that $\int_0^\infty e^{(\lambda-\delta)t}\langle \pi, |\partial_x u|^2 \rangle\,dt < \infty$.

Now, we perform induction. For the convenience, we will use $D$ to mean either $\frac{d}{dx}$ or $\frac{\partial}{\partial x}$. Assume that we have proved that for all $m \leq n-1$

$$\langle \pi, |D^m u|^2 \rangle \leq C_m \exp(-\gamma_m t) \tag{B.3}$$

and that for all $m \leq n$

$$\int_0^\infty e^{\tilde{\lambda}_m t}\langle \pi, |D^m u|^2 \rangle\,dt < \infty. \tag{B.4}$$

We show (B.3)-(B.4) hold for $m \leq n$ and $m \leq n+1$ respectively. Taking the $n$-th order derivative, we have

$$\partial_t D^n u = \mathcal{L}D^n u + g_{n,0}(x)D^{n+1}u + g_{n,1}(x)D^n u + \sum_{m \leq n-1} g_{n,n-m+1}D^m u,$$

where $g_{n,m}(x)$ are smooth functions involving $b, \sigma$ and their derivatives. Multiplying $2\pi D^n u$ and taking integral, we have

$$\partial_t \langle \pi, |D^n u|^2 \rangle \leq -\int_{\mathbb{T}} \Lambda|D^{n+1}u|^2\pi\,dx + C\int_{\mathbb{T}} D^{n+1}u D^n u\pi\,dx$$
$$+ C\langle \pi, |D^n u|^2 \rangle + \sum_{m \leq n-1} C_m \langle \pi, D^m u D^n u \rangle. \tag{B.5}$$

Since $\int_{\overline{\mathbb{T}}} D^{n+1}u D^n u \pi_x \, dx \le \nu \langle \pi, |D^{n+1}u|^2 \rangle + \frac{1}{4\nu}\langle \pi, |D^n u|^2 \rangle$, the $D^{n+1}u$ term is controlled by the first term on the right hand side. The remaining proof is similar as that for $n = 1$, and we omit.

Now that (B.3)-(B.4) hold for all $m \ge 0$. Since $\pi$ is bounded from below, we find that $\|u - \langle \pi, \varphi \rangle\|_{H^k(\mathbb{T})}^2 \le C_n \exp(-\gamma_n t)$. The claims for the decay of $\|u - \langle \pi, \varphi \rangle\|_{C^k}$ follow from Sobolev embedding.

Since $p(x,t) = q(x,t)\pi(x)$ where $q$ satisfies the backward equation for the modified SDE (6.3). The first part of this lemma says that $\|q(\cdot,t) - 1\|_{C^k} \le C \exp(-\gamma_k t)$. Since $\pi$ is smooth on $\mathbb{T}$, we then have $\|\rho(\cdot,t) - \pi\|_{C^k} = \|\pi(q(\cdot,t)-1)\|_{C^k}$ decays to zero exponentially fast. $\qquad\square$

# References

[1] R. J. LeVeque. *Finite difference methods for ordinary and partial differential equations: steady-state and time-dependent problems*, volume 98. SIAM, 2007.

[2] A. Harten, P. D. Lax, and B. v. Leer. On upstream differencing and Godunov-type schemes for hyperbolic conservation laws. *SIAM Review*, 25(1):35–61, 1983.

[3] Mi. G. Crandall and A. Majda. Monotone difference approximations for scalar conservation laws. *Mathematics of Computation*, 34(149):1–21, 1980.

[4] A. Harten and P. D Lax. A random choice finite difference scheme for hyperbolic conservation laws. *SIAM J. Numer. Anal.*, 18(2):289–315, 1981.

[5] A. Harten. On a class of high resolution total-variation-stable finite-difference schemes. *SIAM J. Numer. Anal.*, 21(1):1–23, 1984.

[6] C.-W. Shu. Total-variation-diminishing time discretizations. *SIAM Journal on Scientific and Statistical Computing*, 9(6):1073–1084, 1988.

[7] S. N. Kružkov. First order quasilinear equations in several independent variables. *Mathematics of the USSR-Sbornik*, 10(2):217, 1970.

[8] D. Serre. *Systems of Conservation Laws 1: Hyperbolicity, entropies, shock waves*. Cambridge University Press, 1999.

[9] C. Wang and J.-G. Liu. Positivity property of second-order flux-splitting schemes for the compressible Euler equations. *Discrete Cont. Dyn. Syst. Ser. B*, 3(2):201–228, 2003.

[10] T. M. Liggett. *Continuous Time Markov Processes: An Introduction*, volume 113. American Mathematical Soc., 2010.

[11] L. Miclo. An example of application of discrete Hardy's inequalities. *Markov Process. Related Fields*, 5(3):319–330, 1999.

[12] F. Filbet and M. Herda. A finite volume scheme for boundary-driven convection–diffusion equations with relative entropy structure. *Numerische Mathematik*, 137(3):535–577, 2017.

[13] M. D. Donsker. An invariance principle for certain probability limit theorems. *Mem. Am. Math. Soc.*, 6th, 1951.

[14] M. D. Donsker. Justification and extension of Doob's heuristic approach to the Kolmogorov-Smirnov theorems. *The Annals of Mathematical Statistics*, 23(2):277–281, 1952.

[15] A. Gerardi, F. Marchetti, and A. M. Rosa. Simulation of diffusions with boundary conditions. *Systems & Control Letters*, 4:253–261, 1984.

[16] Harold J. Kushner and Paul Dupuis. *Numerical methods for stochastic control problems in continuous time*, volume 24 of *Applications of Mathematics (New York)*. Springer-Verlag, New York, second edition, 2001.

[17] B. Øksendal. *Stochastic Differential Equations*. Springer, 6 edition, 2003.

[18] X. Mao. *Stochastic Differential Equations and Applications*. Horwood, Chichester, UK, 1997.

[19] R. Khasminskii. *Stochastic Stability of Differential Equations*, volume 66. Springer Science & Business Media, 2011.

[20] S. Meyn and R. L. Tweedie. *Markov Chains and Stochastic Stability*. Cambridge University Press, Cambridge, second edition, 2009.

[21] J. C. Mattingly, A. M. Stuart, and D. J. Higham. Ergodicity for SDEs and approximations: locally Lipschitz vector fields and degenerate noise. *Stoch. Proc. Appl.*, 101(2):185–232, 2002.

[22] A. Eberle. Reflection coupling and Wasserstein contractivity without convexity. *Comptes Rendus Mathematique*, 349(19-20):1101–1104, 2011.

[23] A. Eberle. Reflection couplings and contraction rates for diffusions. *Probab. Theory Relat. Fields*, 166(3-4):851–886, 2016.

[24] L. Gross. Existence and uniqueness of physical ground states. *Journal of Functional Analysis*, 10(1):52–109, 1972.

[25] S. Aida. Uniform positivity improving property, Sobolev inequalities, and spectral gaps. *Journal of Functional Analysis*, 158(1):152–185, 1998.

[26] M. Hino. Exponential decay of positivity preserving semigroups on $L^p$. *Osaka Journal of Mathematics*, 37(3):603–624, 2000.

[27] F. Gong and L. Wu. Spectral gap of positive operators and applications. *Journal de Mathématiques Pures et Appliquées*, 85(2):151–191, 2006.

[28] L. Miclo. On hyperboundedness and spectrum of Markov operators. *Inventiones Mathematicae*, 200(1):311–343, 2015.

[29] P. A. Markowich and C. Villani. On the trend to equilibrium for the Fokker-Planck equation: an interplay between physics and functional analysis. *Mat. Contemp*, 19:1–29, 2000.

[30] P. J. Rabier and C. A. Stuart. Exponential decay of the solutions of quasilinear second-order equations and Pohozaev identities. *Journal of Differential Equations*, 165:199–234, 2000.

[31] L. A. Bagirov. Elliptic equations in unbounded domains. *Sbornik: Mathematics*, 15(1):121–140, 1971.

[32] H Berestycki, L Caffarelli, and L Nirenberg. Further qualitative properties for elliptic equations in unbounded domains. Dedicated to Ennio De Giorgi. *Ann. Scuola Norm. Sup. Pisa Cl. Sci*, 25(1-2):69–94, 1997.

[33] D. Talay. Second-order discretization schemes of stochastic differential systems for the computation of the invariant law. *Stochastics: An International Journal of Probability and Stochastic Processes*, 29(1):13–36, 1990.

[34] L. Wu. Uniformly integrable operators and large deviations for Markov processes. *Journal of Functional Analysis*, 172(2):301–376, 2000.

[35] F. Wang. *Functional Inequalities, Markov Semigroups and Spectral Theory*. Science Press, Beijing, 2005.

[36] L. Gross. Logarithmic Sobolev inequalities. *American Journal of Mathematics*, 97(4):1061–1083, 1975.

[37] J. C. Mattingly, A. M. Stuart, and M. V. Tretyakov. Convergence of numerical time-averaging and stationary measures via Poisson equations. *SIAM J. Numer. Anal.*, 48(2):552–577, 2010.

[38] A. Debussche and E. Faou. Weak backward error analysis for SDEs. *SIAM J. Numer. Anal.*, 50(3):1735–1752, 2012.